\def\draft{n}
\documentclass{amsart}
\usepackage[headings]{fullpage}
\usepackage{amssymb,epic,eepic,epsfig}
\usepackage{graphicx}
\usepackage{texdraw}
\usepackage{url}
\usepackage{amsbsy,amsmath,amsthm,amscd}
\usepackage[bookmarks=true,%
    colorlinks=true,%
    linkcolor=blue,%
    citecolor=blue,%
    filecolor=blue,%
    menucolor=blue,%
    urlcolor=blue,%
    breaklinks=true]{hyperref}

\newtheorem{theorem}{Theorem}[section]
\newtheorem{proposition}{Proposition}[section]
\theoremstyle{definition}

\newtheorem{remark}[proposition]{Remark}

\newtheorem{problem}[proposition]{Problem}

\newcommand{\be}{\begin{equation}}
\newcommand{\ee}{\end{equation}}
\newcommand{\ba}{\begin{aligned}}
\newcommand{\ea}{\end{aligned}}

\def\printname#1{
        \if\draft y
                \smash{\makebox[0pt]{\hspace{-0.5in}
                        \raisebox{8pt}{\tt\tiny #1}}}
        \fi
}

\newcommand{\psdraw}[2]
         {\begin{array}{c} \hspace{-1.3mm}
        \raisebox{-4pt}{\epsfig{figure=draws/#1.eps,width=#2}}
        \hspace{-1.9mm}\end{array}}

\newlength{\standardunitlength}
\catcode`\@=11
\long\def\@makecaption#1#2{%
     \vskip 10pt

\setbox\@tempboxa\hbox{
       \small\sf{\bfcaptionfont #1. }\ignorespaces #2}%
     \ifdim \wd\@tempboxa >\captionwidth {%
         \rightskip=\@captionmargin\leftskip=\@captionmargin
         \unhbox\@tempboxa\par}%
       \else
         \hbox to\hsize{\hfil\box\@tempboxa\hfil}%
     \fi}
\font\bfcaptionfont=cmssbx10 scaled \magstephalf
\newdimen\@captionmargin\@captionmargin=2\parindent
\newdimen\captionwidth\captionwidth=\hsize
\catcode`\@=12

\def\lbl#1{\label{#1}\printname{#1}}

\def\BC{\mathbb C}

\def\ga{\gamma}

\newcommand{\CB}{\mathcal{B}}

\newcommand{\CO}{\mathcal{O}}

\def\IC{{\mathbb C}}

\newcommand{\re}{{\rm e}}
\newcommand{\ri}{{\rm i}}
\newcommand{\rd}{{\rm d}}

\renewcommand{\be}{\begin{equation}}
\renewcommand{\ee}{\end{equation}}
\renewcommand{\ba}{\begin{aligned}}
\renewcommand{\ea}{\end{aligned}}
\newcommand{\ben}{\begin{eqnarray}\displaystyle}
\newcommand{\een}{\end{eqnarray}}

\newcommand{\sectiono}[1]{\section{#1}\setcounter{equation}{0}}

\newdimen\tableauside\tableauside=1.0ex
\newdimen\tableaurule\tableaurule=0.4pt
\newdimen\tableaustep
\def\phantomhrule#1{\hbox{\vbox to0pt{\hrule height\tableaurule width#1\vss}}}
\def\phantomvrule#1{\vbox{\hbox to0pt{\vrule width\tableaurule height#1\hss}}}
\def\sqr{\vbox{%
  \phantomhrule\tableaustep
  \hbox{\phantomvrule\tableaustep\kern\tableaustep\phantomvrule\tableaustep}%
  \hbox{\vbox{\phantomhrule\tableauside}\kern-\tableaurule}}}
\def\squares#1{\hbox{\count0=#1\noindent\loop\sqr
  \advance\count0 by-1 \ifnum\count0>0\repeat}}
\def\tableau#1{\vcenter{\offinterlineskip
  \tableaustep=\tableauside\advance\tableaustep by-\tableaurule
  \kern\normallineskip\hbox
    {\kern\normallineskip\vbox
      {\gettableau#1 0 }%
     \kern\normallineskip\kern\tableaurule}%
  \kern\normallineskip\kern\tableaurule}}
\def\gettableau#1{\ifnum#1=0\let\next=\null\else
\squares{#1}\let\next=\gettableau\fi\next}

\newcommand{\figref}[1]{Fig.~\protect\ref{#1}}

\begin{document}


\title[Asymptotics of the instantons of Painlev\'e I]{
Asymptotics of the instantons of Painlev\'e I}
\author{Stavros Garoufalidis}
\address{School of Mathematics \\
         Georgia Institute of Technology \\
         Atlanta, GA 30332-0160, USA \newline 
         {\tt \url{http://www.math.gatech.edu/~stavros}}}
\email{stavros@math.gatech.edu}
\author{Alexander Its}
\address{Department of Mathematical Sciences \\
402 N. Blackford Street, LD270 \\
Indiana University-Purdue University \\
Indianapolis, IN 46202, USA \newline
         {\tt \url{http://www.science.iupui.edu}}}
\email{itsa@math.iupui.edu}
\author{Andrei Kapaev}
\address{SISSA, via Bonomea, 265, 34136 Trieste, Italy}
\email{kapaev55@mail.ru}
\author{Marcos Mari\~no}
\address{Section de Math\'ematiques et D\'epartement de Physique Th\'eorique \\
         Universit\'e de Gen\`eve \\
         CH-1211 Gen\`eve 4, Switzerland \newline
         {\tt \url{http://www.unige.ch/math/people/marino_en.html}}}
\email{marcos.marino@unige.ch}

\thanks{S.G. and A.I. were supported in part by by NSF. M.M. was supported in part by the Fonds National Suisse.
\newline
1991 {\em Mathematics Classification.} Primary 57N10. Secondary 57M25.
}

\date{\today  }


\begin{abstract}
The $0$-instanton solution of Painlev\'e I is a sequence $(u_{n,0})$
of complex numbers which appears universally in many enumerative 
problems in algebraic geometry, graph theory, matrix models and 2-dimensional 
quantum gravity. The asymptotics of the $0$-instanton $(u_{n,0})$
for large $n$ were obtained by the third author using the Riemann-Hilbert approach. For $k=0,1,2,\dots$,
the $k$-instanton solution of Painlev\'e I 
is a doubly-indexed sequence $(u_{n,k})$ of complex numbers that satisfies an
explicit 
quadratic non-linear recursion relation. The goal of the paper is three-fold:
(a) to compute the asymptotics of the $1$-instanton sequence $(u_{n,1})$ to all orders in $1/n$ by 
using the Riemann-Hilbert method, (b) to present formulas for
the asymptotics of $(u_{n,k})$ for fixed $k$ and to all orders in $1/n$ using resurgent 
analysis, and (c) to confirm numerically the predictions of resurgent analysis. We point out that 
the instanton solutions display a new type of Stokes behavior, induced from the tritronqu\'ee Painlev\'e transcendents, and 
which we call the induced Stokes phenomenon. 
The asymptotics of the 
$2$-instanton and beyond exhibits new phenomena not seen in $0$ and 
$1$-instantons, and their enumerative context is at present unknown.
\end{abstract}

\maketitle

\tableofcontents

\sectiono{Introduction} 
\label{sec.intro}

\subsection{The Painlev\'e I equation and its $0$-instanton solution}
\label{sub.PI}

The  {\em Painlev\'e I equation}
\be
\label{PI}
-\tfrac{1}{6} u'' + u^2 =z
\ee
is a non-linear differential equation with strong integrability properties
that appears universally in various scaling problems; see e.g. \cite{FIKN}.
The formal power series solution
\begin{equation}
\label{eq.u0}
u_0(z)=z^{1/2 } \sum_{n=0}^\infty u_{n,0} z^{-5n/2} \in z^{1/2} \BC[[z^{-5/2}]]
\end{equation}
and correspondingly the sequence $(u_{n,0})$ is the so-called 
$0$-{\em instanton} solution of Painlev\'e I. Substituting \eqref{eq.u0} into
\eqref{PI} collecting terms, and normalizing by setting $u_{0,0}=1$, implies
that $(u_{n,0})$ satisfies the following quadratic recursion relation
\be
\label{diffeq}
u_{n,0} = \frac{25(n-1)^2-1}{48} u_{n-1,0} - {1\over 2} \sum_{\ell=1}^{n-1} 
u_{\ell,0} u_{n-\ell,0}, \qquad u_{0,0} = 1.
\ee
The $0$-instanton solution $(u_{n,0})$ of  Painlev\'e I is a sequence which plays a crucial role in many enumerative 
problems in algebraic geometry, graph theory, matrix models and 2-dimensional 
quantum gravity; see for example \cite{DGZ,GLM,GM}. The leading asymptotics of the 
$0$-instanton 
sequence $(u_{n,0})$ for large $n$ was obtained by the third author using
the {\em Riemann-Hilbert approach}; see \cite{Ka}. In \cite{JK} (see 
also \cite[App.A]{GLM}), asymptotics to all orders in $n$ were
obtained as follows
\be
\label{aas}
u_{n,0} 
\sim A^{-2n+{1\over 2}}  \Gamma\Bigl(2n-{1\over 2} \Bigr)\, 
{S_1 \over \pi \ri} \biggl\{1 + \sum_{l=1}^{\infty} {u_{l,1} A^{l} 
\over \prod_{k=1}^{l} (2n-1/2 -k)} \biggr\}, \qquad n\rightarrow \infty. 
\ee
In this expression, $u_{l,1}$ are the coefficients of the $1$-instanton 
series (defined below), $A$ is the {\it instanton action} 
\be
\label{eq.A}
A= \frac{8\sqrt{3}}{5}
\ee
and
\be
S_1 = -\ri {3^{1\over 4} \over 2{\sqrt{\pi}}}.
\ee
is a Stokes constant. 

\subsection{The Painlev\'e I equation and its $k$-instanton solution}
\label{sub.PIk}

In this paper we study the asymptotics of the $k$-instanton solution
of $(u_{n,k})$ the Painlev\'e I equation. The doubly indexed sequence 
$(u_{n,k})$ for $n,k=0,1,2,\dots$ can be defined concretely by the following
quadratic recursion relation: for $k=1$ we have
  \be
  \label{oneirec}
  u_{n,1}=\frac{8}{25 A n}\Bigl\{ 12\sum_{l=0}^{n-2} u_{l,1} u_{(n+1-l)/2,0} 
-\frac{25}{64}(2n-1)^2 u_{n-1,1} \Bigr\}, \qquad u_{0,1}=1
 . 
 \ee
 while for general $k \ge 2$ we have
 \be
\label{recunl}
 \ba
u_{n,k}&=\frac{1}{12( k^2-1)} \Bigl\{ 12\sum_{l=0}^{n-3} u_{l,k} 
u_{(n-l)/2,0} + 
  6\sum_{m=1}^{k-1} \sum_{l=0}^{n} u_{l,m} u_{n-l,k-m} \\
  & -\frac{25}{64}(2n+k-4)^2 u_{n-2,k}  
  -\frac{25}{16} A k (k+2n-3) 
u_{n-1,k} 
 \Bigr\}.
 \ea
 \ee
It is understood here that $u_{n/2,0}=0$ if $n$ is not an even integer.  

The above recursion defines $u_{n,k}$ in terms of previous $u_{n',k'}$
for $n' <n$ or $n'=n$ and $k' <k$.
The paper is concerned with the asymptotics of the sequence $(u_{n,k})$
for fixed $k$ and large $n$. As we shall see, resurgence analysis
predicts that the asymptotics of $(u_{n,k})$ to all orders in $1/n$ 
is given in terms of three known sequences $(u_{l,k \pm 1})$, 
$(\mu_{l,k \pm 1})$, $(\phi_{l,k-1})$ and two Stokes constants 
$S_1$ and $S_{-1}$; see Equation \eqref{asunkpres}. 
The first constant $S_1$ is known from the 
asymptotics of the $0$-instanton $(u_{n,0})$ and the second one appears 
in the asymptotics of $(u_{n,k})$ for $k \geq 2$ and its exact value is unknown
at present. Two new features appear in the asymptotics of $(u_{n,k})$ for 
$k\ge 2$: 
the constant $S_{-1}$ and the presence of $\log n$ terms. These features 
are absent in the asymptotics of $(u_{n,0})$ and $(u_{n,1})$.

Equation \eqref{recunl} defines but does not motivate the $k$-instanton
solution $(u_{n,k})$ to the Painlev\'e I equation. The motivation comes
from the so-called {\em trans-series solution} of the Painlev\'e I equation.
Trans-series were introduced and studied by \'Ecalle in the eighties;
\cite{Ec1,Ec2,Ec3}. 
The trans-series solution $u(z,C)$ of the Painlev\'e I equation 
is a formal power series in two variables $z$ and $C$
that is defined as follows. Substitute the following expression
\be
\label{uzC}
u(z,C)= \sum_{k\ge0} C^{k} u_{k}(z) 
\ee
into \eqref{PI} and collect the coefficients of $C^k$. It follows that 
$u_k(z)$ satisfy the 
following hierarchy of differential equations, which are non-linear for
$k=0$, linear homogeneous for $k=1$ and linear inhomogeneous for 
$k \geq 2$:
\be
\label{coupledC}
\ba
-\tfrac{1}{6} u_0''+u_0^2 &= z \\
-\tfrac{1}{6} u_k''+\sum_{k'=0}^k u_{k'}u_{k-k'} &= 0, \qquad
k \geq 1.
\ea
\ee
$u_k(z)$ is known as the $k$-instanton solution of \eqref{PI}, and it
has the following structure 
\be
\label{utrans}
u_k (z) =  z^{1/2}\re^{- k A z^{5/4}} \phi_{k}(z), 
\ee
where $A$ is given in (\ref{eq.A}) and
\be\label{phiins}
\phi_{k}(z) =z^{ -5 k /8} \sum_{n\ge 0} u_{n,k} z^{-5n/ 4}.
\ee
Since $C$ is arbitrary we can normalize
\be
  u_{0,1}=1.
  \ee
This motivates our definition of $(u_{n,k})$. The trans-series (\ref{uzC}) 
is what is called a {\it proper} trans-series, in the sense 
that all the exponentials appearing in it are small when 
$z\rightarrow \infty$ along the direction ${\rm Arg}(z)=0$. Therefore, 
the instanton solutions 
$u_k(z)$ give exponentially small corrections to the asymptotic expansion 
(\ref{eq.u0}) and they can be used to construct actual solutions of the 
Painlev\'e I equation in certain sectors; see \cite{C}.  

The instanton solutions of Painlev\'e I also have an important physical 
interpretation in the context of two-dimensional quantum gravity and 
non-critical string theory. As it is well-known (see for 
example \cite{DGZ} and references therein), the Painlev\'e I equation 
appears in the so-called double-scaling limit of random matrices with a 
polynomial potential, and 
it is interpreted as the equation governing the specific heat of a 
non-critical string theory. The asymptotic expansion of the 0-instanton 
solution (\ref{eq.u0}) is nothing but the 
genus or perturbative expansion of this string theory, and $z^{-5/4}$ is 
interpreted as the string coupling constant. The instanton solutions to 
Painlev\'e I correspond to non-perturbative corrections to this expansion. In the context of 
matrix models, they can be interpreted as the double-scaling limit of matrix model instantons, 
which are obtained by eigenvalue tunneling \cite{david,msw}. In the context of non-critical string theory, 
they are due to a special type of D-branes called ZZ branes \cite{martinec,akk}. The asymptotic behavior of the $k$-instanton solution 
is then important in order to understand the full non-perturbative structure of these theories.

\subsection{The predictions of resurgence for the $k$-instanton asymptotics}
\label{sub.predictres}

Our paper consists of three parts.

\begin{itemize}
\item[(a)] A proof of the all-orders asymptotics of the 
$0$ and $1$-instantons of Painlev\'e I, using the Riemann-Hilbert approach. 
\item[(b)] A resurgence analysis of the $k$-instantons of Painlev\'e I
and their asymptotics to all orders.
\item[(c)] A numerical confirmation of the predictions of the resurgence 
analysis.
\end{itemize}
In this section we state the predictions of resurgence analysis for
the asymptotics of $(u_{n,k})$ for fixed $k$ and large $n$. Recall from
Equation \eqref{aas} that the $1/n^k$ asymptotics of $(u_{n,0})$ involve the
1-instanton $u_{\ell,1}$ for $\ell \leq k$, and a Stokes constant $S_1$.
Likewise, the asymptotics of $(u_{n,k})$ for fixed $k\ge1$ and large $n$ 
involves two auxilliary doubly-indexed sequences $(\mu_{n,k})$ and
$(\nu_{n,k})$ and two Stokes constants $S_1$ and $S_{-1}$. 
The sequence $(\nu_{n,k})$ is proportional to the instanton sequence 
 \be
 \label{eq.phirec}
  \nu_{n,k} =\begin{cases} {16\over 5 A} k u_{n,k} &  \qquad k >0 \\
                           0 & \qquad  k=0.
\end{cases}
  \ee
 The sequence $(\mu_{n,k})$ is defined by
\begin{eqnarray*}
\mu_{n,1} &=& (-1)^n u_{n,1} \\ 
\mu_{2n,2} &=&
-\sum_{l=0}^{n-2} \mu_{2l,2} u_{n-l,0} -\sum_{l=0}^{2n}  
(-1)^l u_{l,1} u_{ 2n-l,1} +{25\over 192}(2n-1)^2 \mu_{2(n-1),2} \\
\mu_{2n+1,2} &=& 0 \\
\mu_{n,3}&=&{8 \over 25 A (n+1)}\biggl\{ -{25\over 64} (2n+1)^2 \mu_{n-1,3} 
+ 12 \sum_{l=0}^{n-2} u_{(n+1-l)/2,0}\mu_{l,3} \\
& & \qquad + 12 \sum_{m=1}^{2} \sum_{l=0}^{n+1} u_{n+1-l,m}\mu_{l,3-m} 
+ {25\over 16} (2n+1) \nu_{n,1} + {25A\over 8} \nu_{n+1,1}\biggr\} 
\end{eqnarray*}
for $k=1,2,3$, and by
\begin{eqnarray*}
\mu_{n,k}&=& \frac{1}{12( k-1)(k-3)} \Bigl\{ 12\sum_{l=0}^{n-3} \mu_{l,k} 
u_{(n-l)/2,0} + 
  12 \sum_{m=1}^{k-1} \sum_{l=0}^{n} \mu_{l,m} u_{n-l,k-m} 
-\frac{25}{64}(2n+k-4)^2 \mu_{n-2,k}   \\
  & &
  -\frac{25}{16} A k (k+2n-3) 
\mu_{n-1,k} + {25\over 16} (k+2n-4) \nu_{n-1,k-2} 
+ {25A\over 8}(k-2) \nu_{n,k-2}
 \Bigr\}
\end{eqnarray*}
for $k \geq 4$. We apologize for the lengthy formulas that define
the doubly-indexed sequences $(\mu_{n,k})$ and $(\nu_{n,k})$. 
As in Section \ref{sub.PIk}, there is a simple resurgence analysis
explanation of these sequences given in detail in Section \ref{sec.trans}.
An analysis based on trans-series solutions and resurgent properties suggests the following asymptotic result for the 
coefficients $u_{n,k}$, for fixed $k\ge1$, in the limit $n \rightarrow \infty$:
\be
\label{asunkpres}
\ba
u_{n,k} & \sim_n  A^{-n+1/2} {S_1\over 2 \pi\ri}  
\Gamma\bigl(n-1/2 \bigr)\, 
\biggl\{ (k+1) u_{0,k+1} +(-1)^n \mu_{0,k+1}  
+ \sum_{l=1}^{\infty} {((k+1)u_{l,k+1} +(-1)^{n+l} 
\mu_{l,k+1}) A^{l} \over \prod_{m=1}^{l} 
(n-1/2 -m)} \biggr\} \\ 
& + 
(-1)^n (k-1) A^{-n-2/3} {S_{-1} \over 2\pi\ri}  
\Gamma\bigl(n+1/2 \bigr)\, 
\biggl\{ u_{0,k-1}+ \sum_{l=1}^{\infty} {u_{l,k-1} (-A)^{l} \over 
\prod_{m=1}^{l}
(n+1/2-m)} \biggr\}\\
&- (-1)^n A^{-n-1/2} {S_{1} \over 2\pi\ri}  
\Gamma\bigl(n+1/2 \bigr)\, (\log n -\log A) \biggl\{ \nu_{0,k-1}
+ \sum_{l=1}^{\infty} {\nu_{l,k-1} (-A)^{l} \over \prod_{m=1}^{l}
(n+1/2-m)} \biggr\}\\
 &-(-1)^n A^{-n-1/2} {S_{1} \over 2\pi\ri}  \Gamma\bigl(n+1/2 
\bigr)\,  \biggl\{ (\psi(n+1/2) -\log n) \nu_{0,k-1} 
+ \sum_{l=1}^{\infty} {\psi(n+1/2-l) 
-\log n\over \prod_{m=1}^{l}
(n+1/2-m)} \nu_{l,k-1} (-A)^{l} \biggr\},
\ea
\ee
where 
$$
\psi(z)=\frac{\Gamma'(z)}{\Gamma(z)}
$$
is the logarithmic derivative of the $\Gamma$ function; see \cite{O}.
In the above Equation, the convention is that
$u_{n,k}=\nu_{n,k}=0$ for $k =0$. 

The above formula gives an asymptotic expansion for $u_{n,k}$ for $n$ large, involving the coefficients $u_{l,k\pm 1}$, $\mu_{l, k\pm1}$ and 
$\nu_{l, k\pm 1}$. Up to the overall factor $A^{-n}\Gamma (n+1/2)$, this expansion involves terms of the form $1/n^l$ and $\log n/n^l$. In order to get the $m$-th first terms of this asymptotic expansion for $u_{n,k}$, we only need to know the coefficients $u_{l,k\pm 1}$, $\mu_{l, k\pm1}$, $\nu_{l, k\pm 1}$ with $l$ up to $m$. In particular, (\ref{asunkpres}) gives an efficient method to obtain the asymptotic behavior of the $u_{n,k}$ for fixed $k$ at large $n$. 

As a concrete and important example which also clarifies the above remarks, let us look at $k=1$. In this case, only the first line of Equation \eqref{asunkpres}
contributes and resurgence analysis predicts  that for even $n$, we have
\be
\label{ueven1}
u_{2n,1} \sim A^{-2n+1/2} {S_1\over 2 \pi\ri}  \Gamma\bigl(2n -1/2 
\bigr)\biggl\{ 2u_{0,2} +  \mu_{0,2}  +  \sum_{l=1}^{\infty} 
{(2 u_{l,2} +(-1)^l\mu_{l,2}) A^{l} \over \prod_{m=1}^{l} 
(2n-1/2 -m)}\biggr\} , \qquad n \rightarrow \infty, 
\ee
while for odd $n$ we have
\be
\label{uodd1}
u_{2n+1,1} \sim A^{-2n-1/2} {S_1\over 2 \pi\ri}  
\Gamma\bigl(2n +1/2 \bigr)\biggl\{ 2u_{0,2} -  \mu_{0,2}  
+  \sum_{l=1}^{\infty} {(2 u_{l,2} -(-1)^l\mu_{l,2}) \lambda^{l} \over 
\prod_{m=1}^{l} 
(2n+1/2 -m)}\biggr\}, \qquad n \rightarrow \infty.
\ee
This prediction will be proved in section 3 by using the Riemann--Hilbert approach. The coefficients $u_{0,2}$ and $\mu_{0,2}$ can be easily calculated from the 
recursion (\ref{recunl}):
\be
\label{exonei}
u_{0,2}={1\over 6}, \qquad \mu_{0,2}=-1.
\ee
In fact one can obtain closed formulae for $u_{0,k}$ and $\mu_{0,k}$ for all $k$ (see (\ref{onel}) and (\ref{mugen}), respectively). The explicit result (\ref{exonei}) gives the leading 
asymptotics, 
\be
u_{n,1}  \sim \left( {1\over 3}  - (-1)^n \right) A^{-n+1/2} {S_1\over 2 \pi\ri}  \Gamma\bigl(n-1/2 
\bigr), \qquad n\rightarrow \infty.. 
\ee
Concrete examples of the implications of (\ref{asunkpres}) for the asymptotics of the $u_{n,k}$ 
are given in section 5, where the formula is tested numerically. 


 For some partial results on the Painlev\'e I equation and the asymptotics of the coefficients of 
 the 0-instanton solution, see also \cite{C,CC,CK2}. To the best of our knowledge, a rigorous computation
of the Stokes constant $S_1$ has only been achieved via the Riemann-Hilbert
approach \cite{Ka} or its earlier version - the isomonodromy method \cite{Ka1,Tak}. 
If at all possible, a computation of the constants $S_1$ and $S_{-1}$ using
resurgence analysis would be very interesting.

\subsection{The asymptotics of the 1-instanton solution via the
Riemann-Hilbert method}
\label{sub.1inst}

Recall from Equations \eqref{coupledC} and \eqref{utrans} 
that the generating series $u_1(z)$ of the 1-instanton $(u_{n,1})$
solution to Painlev\'e I satisfies the differential equation
\be
\label{one000}
u_1''(z)=12 u_1(z) u_0(z)
\ee
where $u_0(z)$ is the generating series of the 0-instanton solution
$(u_{n,0})$ solution to Painlev\'e I. There are five explicit 
{\it tritronqu\'ees} meromorphic solutions of Painlev\'e I equation
asymptotic to $u_0(z)$ in appropriate sectors in the complex plane,
discussed in Section \ref{reviewtri}.
To study the asymptotics of the 1-instanton $(u_{n,1})$, consider the 
linear homogenous differential equation
\begin{equation}
\label{eq.gammauv}
v''=\gamma u v
\end{equation}
where $\gamma$ is a constant and $u$ is a tritronqu\'ee solution to 
Painlev\'e I. It is easy to see that Equation \eqref{eq.gammauv} has
two linearly independent formal power series solutions $v_f^{\pm}$
of the form 
\begin{equation}
\label{vformal}
v_f^{\pm}(z)=z^{-1/8}e^{\pm 4 \frac{\sqrt{\ga}}{5} z^{5/4}}
\sum_{n=0}^{\infty}b_n^{\pm}(\ga)z^{-5n/4}
\end{equation}
where $b_n^{\pm}=b_n^{\pm}(\ga)$ is given by
\begin{equation}
b_n^{\pm}=\frac{1}{2n}
\biggl\{
\frac{5 b_{n-1}^{\pm}}{4 \sqrt{\ga}}\bigl(n-\tfrac{1}{10}\bigr)
\bigl(n-\tfrac{9}{10}\bigr)
\mp \frac{4 \sqrt{\ga}}{5} 
\sum_{m=1}^{[\frac{n+1}{2}]}u_{m,0}b_{n+1-2m}^{\pm}
\biggr\},\qquad b_0^{\pm}=1.
\end{equation}
Note that $b_n^{\pm}(\ga)$ is a polynomial in $\ga^{1/2}$ and 
$b_n^{-}(12)=u_{n,1}$. The next theorem gives the asymptotic expansion
of $b_n^{-}$ when $\ga > 3$. Let

\begin{equation}
\label{eq.B}
B=\frac{4}{5} \sqrt{\ga}.
\end{equation}

\begin{theorem}
\label{thm.bn}
For $n$ large and $\ga > 3$, we have:
\begin{equation}
\label{as.bn}
b_n^-(\ga)=
\frac{4\cdot3^{1/4}}{25\pi^{3/2}}
\frac{A(1+(-1)^n)-2B(1-(-1)^n)}{4B^2-A^2}
\gamma
A^{-n-\frac{1}{2}}
\Gamma(n-\tfrac{1}{2})
\bigl(1+{\mathcal O}(n^{-1})\bigr).
\end{equation}
\end{theorem}
The proof of the theorem uses explicitly the 5 {\it tritronqu\'ees} 
solutions of Painlev\'e I equation and the asymptotic expansion of 
their difference in 5 sectors of the complex plane, 
computed by the Riemann-Hilbert approach to Painlev\'e I. 
An analytic novelty of Theorem \ref{thm.bn} is the rigorous computation
of a Stokes constant which is independent of $\ga$, when $\ga > 3$. 

\begin{remark}
There are two special values of $\gamma$ in Equation \eqref{eq.gammauv}. When
$\gamma=12$ and $v$ is a solution of Equation \eqref{eq.gammauv}, it follows that $v$ is the 1-instanton solution of Painleve I and Theorem \ref{thm.bn} implies Equations  \eqref{ueven1} and \eqref{uodd1}. When $\gamma=3/4$ and $v$ is a solution of \eqref{eq.gammauv}, then $\tilde{v} = -2v'/v$ is a solution of the Riccati equation  
\be
2\tilde v'-\tilde v^2+3 u=0
\ee
studied in \cite{GM}.
\end{remark}

\subsection{Acknowledgements}
M.M. would like to thank Jean--Pierre Eckmann, Albrecht Klemm, Sara Pasquetti, 
Pavel Putrov, Leonardo Rastelli, Marco Rauch, and Ricardo Schiappa for discussions.

\sectiono{Asymptotics of the $0$-instanton of Painlev\'e I}
\label{as.0inst}

In this section we review and extend the calculation of the asymptotics of the $0$-instanton 
sequence $(u_{n,0})$ presented in \cite{Ka}:

\begin{itemize}
\item[(a)] Use the 5 tritronqu\'ee solutions to Painlev\'e I and their
analytic properties obtained by the Riemann-Hilbert method, as an input.
\item[(b)] No two tritronqu\'ee solutions are equal on a sector, however their
difference is exponentially small. Using the 5 tritronqu\'ee solutions and 
their differences, define a piece-wise analytic function in the complex plane
which has (\ref{eq.u0}) as its uniform asymptotics in the neighborhood of 
infinity.
\item[(c)] Apply a mock version of the Cauchy integral formula to obtain
an exact integral formula for $u_{n,0}$ in terms of the jumps of the piece-wise
analytic function defined above. The knowledge of the explicit value
of the relevant Stokes' multiplier (available due to the Riemann-Hilbert 
analysis) yields the large $n$ asymptotics of $u_{n,0}$. In fact, we extend the result of \cite{Ka} and we obtain the asymptotics 
to all orders in $1/n$. 
\end{itemize}
Note that this method consists of working entirely in the $z$-plane 
(and not in the Borel plane), in all sectors simultaneously, and our glued
function is only piece-wise analytic. This method is different from the 
method of Borel transforms analyzed in detail in \cite{CK2,CC}.

\subsection{A review of the tritronqu\'ee solutions of Painlev\'e I}
\lbl{reviewtri}

In this section, to simplify our notation, we will replace $u_{n,0}$ by 
$a_{n}$, and use the symbol $u_n$ {\em not} for the $n$-th term in the 
formal trans-series expansion, but for certain 
exact solutions of the first Painlev\'e equation specified below.

Recall that equation PI \eqref{PI} admits a formal 0-parameter solution
with the power expansion $u_f(z)$  (cf. (\ref{eq.u0})),
\begin{equation}
\label{u0_formal}
u_f(z)=z^{1/2}\sum_{n=0}^{\infty}a_n z^{-5n/2},\quad
a_0=1,\quad
a_1=-1/48,\quad
a_{n+1}=\frac{25n^2-1}{48}a_n-\frac{1}{2}\sum_{m=1}^na_ma_{n+1-m}.
\end{equation}

\begin{remark}
\label{u_y_rel}
In \cite{Ka}, the Painlev\'e I was studied in the following form:
\be
\label{PIK}
y_{xx}=6y^2+x.
\ee
The change of variables $u=6^{2/5}y$ and 
$z=e^{-i\pi}6^{-1/5}x$ transforms \eqref{PI} to \eqref{PIK}.
Comparing with \cite{Ka}, we also find more 
convenient to modify the indexes for the tritronqu\'ee solutions.
\end{remark}

\begin{figure}[htpb]
$$
\psdraw{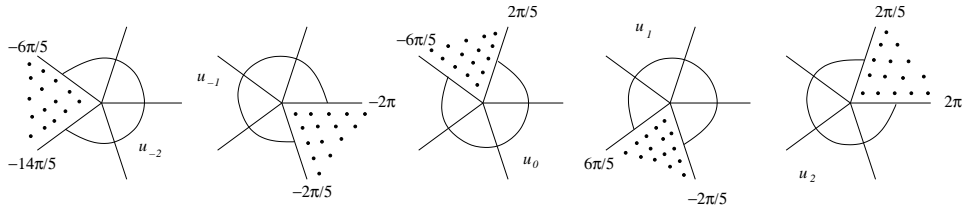}{5in}
$$
\caption{The sectors of the $z$-complex plane where the tritronqu\'ee 
solutions (\ref{u0_gen}) $u_{-2}$, $u_{-1}$, $u_0$, $u_1$, $u_2$ are 
represented by the formal series $u_f$.In the
dotted sectors, the asymptotics at infinity of the tritronqu\'ee solutions 
is elliptic.}
\label{tritronquee}
\end{figure}

Using the Riemann-Hilbert approach (see \cite[Ch.5.2]{FIKN} and also 
\cite{Ka}), it can be shown that there 
exist five different genuine meromorphic solutions of \eqref{PI} with the 
asymptotic power expansion \eqref{u0_formal} in one of the sectors of the 
$z$-complex plane of opening $8\pi/5$, see \figref{tritronquee}:

\begin{eqnarray}
\nonumber
u_0(z) &\sim & u_f(z),\quad z\to\infty,\quad
\arg z\in\bigl(-\tfrac{6\pi}{5},\tfrac{2\pi}{5}\bigr),
\\
\label{u0_gen}
u_k(z) &=& e^{-i\frac{8\pi}{5}k}u_0\bigl(e^{-i\frac{4\pi}{5}k}z\bigr)
\sim u_f(z),\quad z\to\infty,\quad
\arg z\in\bigl(-\tfrac{6\pi}{5}+\tfrac{4\pi}{5}k,\tfrac{2\pi}{5}
+\tfrac{4\pi}{5}k\bigr),
\\
u_{k+5}(z) &=& u_k(z),
\end{eqnarray}
where in $u_f(z)$  appearing in (\ref{u0_gen}) the branches of $z^{1/2}$ and $z^{-5n/2}$ 
are defined according to the rule,
\begin{equation}\label{rule}
z^{1/2} =  \sqrt{|z|}e^{\arg z/2}, \quad
 \arg z \in\bigl(-\tfrac{6\pi}{5}+\tfrac{4\pi}{5}k,\tfrac{2\pi}{5}
+\tfrac{4\pi}{5}k\bigr),
\quad z^{-5n/2} = \left(z^{1/2}\right)^{-5n}.
\end{equation}
The asymptotic formula (\ref{u0_gen}) is understood in the usual sense. That is,
for every $\epsilon >0$ and natural number $N$ there exist positive constants
$C_{N,\epsilon}$ and $R_{N,\epsilon}$ (depending on $N$ and $\epsilon$ only)
such that,

\begin{equation}\label{new1}
\Bigl|u_k(z) - z^{1/2}\sum_{n=0}^{N}a_nz^{-5n/2}\Bigr| < C_{N,\epsilon}
|z|^{-\frac{5N+4}{2}},
\end{equation}

$$
\forall z : |z| \geq R_{N,\epsilon}, \quad \arg z \in 
\bigl[-\tfrac{6\pi}{5}+\tfrac{4\pi}{5}k +\epsilon,\tfrac{2\pi}{5}
+\tfrac{4\pi}{5}k -\epsilon\bigr].
$$
Here,  the branches of $z^{1/2}$ and $z^{-5n/2}$ 
are again defined according to (\ref{rule}).
Estimate (\ref{new1}) is proven in \cite{Ka}. 
\begin{remark}\label{poles}
We  note that estimate (\ref{new1}) 
implies, in particular,  that each tritronqu\'ee solution $u_k(z)$ might have in the
sector $ \bigl[-\tfrac{6\pi}{5}+\tfrac{4\pi}{5}k +\epsilon,\tfrac{2\pi}{5}
+\tfrac{4\pi}{5}k -\epsilon\bigr]$ only finitely many  poles which lay inside
the  circle of radius $R_{0, \epsilon}$. We  also note that the number $R_{N,\epsilon}$
in  estimate (\ref{new1}) can be replaced by $R_{0,\epsilon}$ for every $N$.
\end{remark}

Furthermore, in \cite{Ka} it is shown that  the exponential small difference between the tritronqu\'ee 
solutions $u_{k}(z)$ and $u_{k+1}(z)$,
within the common sector where they have 
the identical asymptotics $u_f(z)$ in all orders, 
admits the following explicit asymptotic description:
\begin{gather}
\label{quasi_linear_Stokes}
z\to\infty,\quad
\arg z\in\bigl(-\tfrac{2\pi}{5},\tfrac{2\pi}{5}\bigr)\colon\quad
u_1(z)-u_0(z)=i\frac{3^{1/4}}{2\sqrt\pi}\,z^{-1/8}e^{-\frac{8\sqrt3}{5}z^{5/4}}
\Bigl(1+{\mathcal O}\bigl(z^{-5/4}\bigr)\Bigr)
\\
\nonumber
z\to\infty,\quad
\arg z\in\bigl(-\tfrac{2\pi}{5}+\tfrac{4\pi}{5}k,
\tfrac{2\pi}{5}+\tfrac{4\pi}{5}k\bigr)\colon
\\
\nonumber
u_{k+1}(z)-u_k(z)=e^{i\frac{\pi}{2}(k+1)}
\frac{3^{1/4}}{2\sqrt\pi}\,z^{-1/8}
e^{(-1)^{k+1}\frac{8\sqrt3}{5}z^{5/4}}
\Bigl(1+{\mathcal O}\bigl(z^{-5/4}\bigr)\Bigr).
\end{gather}
Here again the asymptotic relations mean that the differences $u_{k+1}(z)-u_k(z)$ admit the representation,
\begin{equation}\label{new2}
u_{k+1}(z)-u_k(z)=e^{i\frac{\pi}{2}(k+1)}
\frac{3^{1/4}}{2\sqrt\pi}\,z^{-1/8}
e^{(-1)^{k+1}\frac{8\sqrt3}{5}z^{5/4}}
\Bigl(1+r(z)\Bigr),
\end{equation}
whith the error term $r(z)$ satisfying  the estimate (cf. (\ref{new1}),
\begin{equation}\label{new3}
|r(z)| < C_{\epsilon}
|z|^{-5/4},
\end{equation}
$$
\forall z : |z| \geq R_{\epsilon}, \quad \arg z \in 
\bigl[-\tfrac{2\pi}{5}+\tfrac{4\pi}{5}k +\epsilon,\tfrac{2\pi}{5}
+\tfrac{4\pi}{5}k -\epsilon\bigr],
$$
with the positive constants $C_{\epsilon}$ and $R_{\epsilon}$ depending
on $\epsilon$ only. In fact, we can take $ R_{\epsilon} = R_{0,\epsilon}$.

When dealing with the tritronqu\'ee solutions $u_k(x)$, it is convenient 
to assume that  $k$ can take any integer value, simultaneously  remembering that
$u_{k+5}(x) = u_{k}(x)$, see \eqref{u0_gen}. In other words, in the notation
 $u_k(x)$ we shall assume that
 \begin{equation}\label{module}
 k \in {\Bbb Z}, \quad \mbox{mod} \quad 5,
 \end{equation}
unless $k$ is particularly specified, as in \eqref{hat_uN_def} below.
\begin{remark}
\label{tritrong000}
The existence of the tritronqu\'ee solutions $u_k(x)$ can be proved without 
use of the Riemann-Hilbert method (see e.g. \cite{JK} and references therein).
Indeed, these solutions had already been known to Boutroux. 
The Riemann-Hilbert method is needed for the exact evaluation of the pre-exponent
numerical coefficient in the jump-relations (\ref{quasi_linear_Stokes}), i.e.
for an explicit description of the {\em quasi-linear Stokes' phenomenon} 
exhibited by the first Painlev\'e equation.
\end{remark}

\subsection{Gluing the five tritronqu\'ee solutions together}
\lbl{sub.gluing}

The main technical part of the approach of \cite{Ka} to the 0-instanton asymptotics
is a piece-wise meromorphic function with uniform asymptotics \eqref{u0_formal} at 
infinity. This function is constructed from the collection of the 
tritronque\'ee solutions $u_k(z)$, $k=0,\pm1,\pm2$, as follows.

First, observe that the change of independent variable $z=t^2$ turns 
\eqref{u0_formal} 
into the non-branching series,
\begin{equation}
\label{u_t_Laurent}
\hat u_f(t):=u_f(t^2)=\sum_{n=0}^{\infty}a_n t^{-5n+1}.
\end{equation}
Multiplying \eqref{u_t_Laurent} by $t^{5N-2}$, we obtain
the formal series
\begin{equation}
\label{hat_uf_N_def}
\hat u_f^{(N)}(t)\equiv \hat u_f(t)t^{5N-2} = P_{5N-1}(t)
+a_N t^{-1}
+\sum_{n=N+1}^{\infty}a_n t^{-5(n-N)-1},\quad
P_{5N-1}(t)=\sum_{n=1}^{N}a_{N-n}t^{5n-1}.
\end{equation}
Let us perform the similar operation with the solutions $u_{k}(z)$
\begin{gather}
\label{hat_uk_def}
\hat u_k(t)=u_k(t^2),
\\
\label{hat_uj_N_def}
\hat u_k^{(N)}(t)=\hat u_k(t)t^{5N-2}-P_{5N-1}(t).
\end{gather}
Observe that,
$$
\hat u_k^{(N)}(t)=
a_N t^{-1}
+{\mathcal O}_{N}\bigl(t^{-6}\bigr),\quad
t\to\infty,\quad
\arg t\in\bigl(-\tfrac{3\pi}{5}+\tfrac{2\pi}{5}k,\tfrac{\pi}{5}
+\tfrac{2\pi}{5}k\bigr).
$$
The symbol ${\mathcal O}_{N}$ indicates that the relevant positive  constants in the estimate
depend on $N$. 
More precisely,  estimate (\ref{new1}) implies that
\begin{equation}\label{new4}
\Bigl|\hat u_k^{(N)}(t) - a_Nt^{-1}\Bigr| < C_{N,\epsilon} |t|^{-6},
\end{equation}

$$
\forall t : |t| \geq  R^{1/2}_{0,\epsilon}, \quad \arg t \in 
\bigl[-\tfrac{3\pi}{5}+\tfrac{2\pi}{5}k +\epsilon,\tfrac{\pi}{5}
+\tfrac{2\pi}{5}k -\epsilon\bigr].
$$
Moreover, in view of Remark \ref{poles},  we conclude that each function 
$\hat u_k^{(N)}(t)$ might have in the sector $
\bigl[-\tfrac{3\pi}{5}+\tfrac{2\pi}{5}k +\epsilon,\tfrac{\pi}{5}
+\tfrac{2\pi}{5}k -\epsilon\bigr]$ only  finitely many poles whose
number does not depend on $N$, and they all lay inside the 
circle of radius $ R^{1/2}_{0,\epsilon}$. Indeed,  by (\ref{hat_uj_N_def}),    the
possible poles must coincide, for every $N$, with the possible  
poles of the corresponding function $\hat u_k(t)$.

\begin{figure}[htpb]
$$
\psdraw{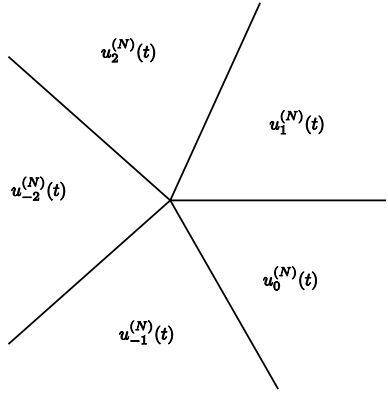}{2in}
$$
\caption{The sectorially meromorphic function $\hat u^{(N)}(t)$ 
\eqref{hat_uN_def}.}
\label{NRH0}
\end{figure}

Put
\begin{equation}
\label{hat_uN_def}
\hat u^{(N)}(t)=\hat u_k^{(N)}(t),\quad
\arg t\in\bigl(-\tfrac{2\pi}{5}+\tfrac{2\pi}{5}k,\tfrac{2\pi}{5}k\bigr),
\quad k = 0, \pm 1, \pm 2.
\end{equation}
These equations determine $\hat u^{(N)}(t)$ as a
sectorially meromorphic function $\hat u^{(N)}(t)$ discontinuous 
across the rays
$\arg t=\frac{2\pi}{5}k$, $k=0,\pm1,\pm2$, see \figref{NRH0}.
Let us choose and then {\it fix} the parameter $\epsilon\equiv \epsilon_0$ in
such a way that the closure of  each sector $\bigl(-\tfrac{2\pi}{5}+\tfrac{2\pi}{5}k,\tfrac{2\pi}{5}k\bigr)$
is included in the corresponding sector $\bigl[-\tfrac{3\pi}{5}+\tfrac{2\pi}{5}k +\epsilon,\tfrac{\pi}{5}
+\tfrac{2\pi}{5}k -\epsilon\bigr]$. Then, the function 
$\hat u^{(N)}(t)$ would have no more than   finitely many poles in the whole complex plane, the
number of poles would be the same for all $N$, and they all lay inside the
circle of radius
\begin{equation}\label{R0def}
R_0 =  R^{1/2}_{0,\epsilon_0}
\end{equation}
In addition, the following 
{\it uniform} asymptotics at infinity takes place,
\begin{equation}
\label{hat_uN_as}
\hat u^{(N)}(t)=
a_N t^{-1}
+{\mathcal O}_{N}\bigl(t^{-6}\bigr),\quad
t\to\infty.
\end{equation}
The exact meaning of this estimate is the existence of
the positive constant $C_{N} (\equiv C_{N,\epsilon_0})$ such that,
\begin{equation}\label{new5}
\Bigl|\hat u^{(N)}(t) - a_Nt^{-1}\Bigr| < C_{N} |t|^{-6},
\end{equation}
$$
\forall t: |t| > R_{0}.
$$
Here, $R_0$ is defined in (\ref{R0def}), and it does not depend on $N$. We remind
also the reader that the circle of radius $R_0$  contains all the possible poles of $\hat u^{(N)}(t)$.

The jumps of the function $\hat u^{(N)}(t)$ across the rays $\arg t=\frac{2\pi}{5}k$, $k=0,\pm1,\pm2$,
oriented towards infinity,
are described by the equations,
\begin{equation}
\label{hat_uN_jumps}
\hat u^{(N)}_+(t)-\hat u^{(N)}_-(t)=\hat U^{(N)}(t),\quad
\arg t=\tfrac{2\pi}{5}k,\quad
k=0,\pm1,\pm2,
\end{equation}
where $\hat u^{(N)}_+(t)$ and $\hat u^{(N)}_-(t)$ are the limits of 
$\hat u^{(N)}(t)$ as we approach the rays from the left and from 
the right, respectively. (We note that $\hat u_6^{(N)}(t) = \hat u_1^{(N)}(t)$,
as it follows from \eqref{u0_gen}).
By virtue of (\ref{quasi_linear_Stokes}), the jump functions
$\hat U^{(N)}(t)$ satisfy the estimate,
\begin{equation} 
\label{hat_UN_as}
\hat U^{(N)}(t)\bigr|_{\arg t=\tfrac{2\pi}{5}k}=
e^{i\frac{\pi}{2}(k+1)}
\frac{3^{1/4}}{2\sqrt\pi}\,t^{5N-9/4}
e^{(-1)^{k+1}\frac{8\sqrt3}{5}t^{5/2}}
\Bigl(1+{\mathcal O}\bigl(t^{-5/2}\bigr)\Bigr),
\end{equation}
which, again, means (cf. (\ref{new2}) and (\ref{new3})) that 
 \begin{equation} 
\label{new6}
\hat U^{(N)}(t)\bigr|_{\arg t=\tfrac{2\pi}{5}k}=
e^{i\frac{\pi}{2}(k+1)}
\frac{3^{1/4}}{2\sqrt\pi}\,t^{5N-9/4}
e^{(-1)^{k+1}\frac{8\sqrt3}{5}t^{5/2}}
\Bigl(1+r(t)\bigr)\Bigr),
\end{equation}
with
$$
|r(t)| < C|t|^{- 5/2}, \quad \forall t : |t| > R_0, \quad \arg t = \tfrac{2\pi}{5}k,
$$
where the constant  $C (\equiv C_{\epsilon_0})$, similar to $R_0$,  is a numerical constant not depending on $N$.

\subsection{An integral formula for the $0$-instanton coefficients and their asymptotic
expansion}
\begin{figure}[htpb]
$$
\psdraw{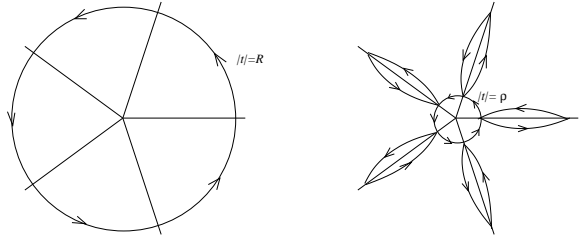}{3in}
$$
\caption{Deformation of the contour of integration for computation of the 
$0$-instanton $N$-large asymptotics}
\label{collapsed0}
\end{figure}

Consider the integral,
$$
\frac{1}{2\pi i}\oint_{|t|=R}\hat u^{(N)}(t)\,dt,
$$
of the function $\hat u^{(N)}$ along the circle of radius $R$ 
centered at the origin and counter-clockwise oriented.  For sufficiently
large $R$ we can apply to the integrand estimate (\ref{new5}). Noticing that
the constant $C_N$ is the same along the whole circle, we conclude that
\begin{equation}
\label{residue0_infinity}
\frac{1}{2\pi i}\oint_{|t|=R}\hat u^{(N)}(t)\,dt=a_N+{\mathcal O}_{N}(R^{-5}),
\end{equation}
that is, 
\begin{equation}\label{new10}
\Bigl|\frac{1}{2\pi i}\oint_{|t|=R}\hat u^{(N)}(t)\,dt-a_N\Bigr| < C_{N}R^{-5}, \quad
\forall R > R_{0}.
\end{equation}
On the other hand, since $\hat u^{(N)}(t)$ can have only a finite number 
of poles lying inside of the circle with the radius $R_0$,  the circular contour of integration can be deformed to the sum of
the circle of smaller radius $|t|=\rho \geq R_0$, still containing inside all the possible poles of 
$\hat u^{(N)}(t)$,
and positive and negative sides of segments of the rays 
$\arg t=\frac{2\pi}{5}k$,
see \figref{collapsed0}.

In other words, taking into account \eqref{hat_uN_jumps}, we have that
\begin{multline}
\label{residue0_integral_app}
a_N=\frac{1}{2\pi i}\oint_{|t|=R\gg1}\hat u^{(N)}(t)\,dt+{\mathcal O}_{N}(R^{-5})=
\\
=-\frac{1}{2\pi i}\sum_{k=-2}^2\int_{\rho e^{i\frac{2\pi}{5}k}}^{R 
e^{i\frac{2\pi}{5}k}}
\hat U^{(N)}(t)\,dt
+\frac{1}{2\pi i}\oint_{|t|=\rho}\hat u^{(N)}(t)\,dt
+{\mathcal O}_{N}(R^{-5})=
\\
=-\frac{1}{2\pi i}\sum_{k=-2}^2\int_{\rho e^{i\frac{2\pi}{5}k}}^{R 
e^{i\frac{2\pi}{5}k}}
\hat U^{(N)}(t)\,dt
+\frac{1}{2\pi i}\oint_{|t|=\rho}\Bigl(\hat u^{(N)}(t) + P_{5N-1}(t)\Bigr) \,dt
+{\mathcal O}_{N}(R^{-5}),
\end{multline}
where we notice  that the integral of the 
polynomial $P_{5N-1}(t)$ can be indeed added to the right hand side since it  is zero.
Let us now notice that from  definition (\ref{hat_uN_def}) of the function $\hat u^{(N)}(t)$
it follows that
\begin{equation}\label{new11}
\hat{u}^{(N)}(t)  + P_{5N-1}(t) = \hat{u}(t)t^{5N-1}, 
\end{equation}
where $\hat{u}(t)$ is defined by the equations (cf. (\ref{hat_uN_def})),
\begin{equation}
\label{hat_u_def}
\hat u(t)=\hat u_k(t),\quad
\arg t\in\bigl(-\tfrac{2\pi}{5}+\tfrac{2\pi}{5}k,\tfrac{2\pi}{5}k\bigr),
\quad k = 0, \pm 1, \pm 2.
\end{equation}
By exactly the same reasons as in the case of estimate (\ref{new5}), we conclude from
(\ref{new1}) that
\begin{equation}\label{new12}
|\hat u(t)| < C^{(1)} |t|,
\quad
\forall t: |t| > R_0,
\end{equation}
with a numerical constant $C^{(1)}$ this time{\footnote{The constant $C^{(1)}$ 
can be taken equal to
$$
C^{(1)} = 1 + \frac{C_{0,\epsilon_0}}{R^{5}_0}.
$$}}  independent of $N$. Estimate (\ref{new12}), together with 
(\ref{new11}) implies the inequality,
$$
|\hat u^{(N)}(t) + P_{5N-1}(t)|_{|t| = \rho} < C^{(1)} \rho^{5N-1},
$$
where  $\rho$ is assumed to be
a fixed positive number satisfying,
$$
R_0\leq\rho < R.
$$
Therefore, formula
(\ref{residue0_integral_app}) can be transformed into the formula,
\begin{equation}\label{new13}
a_N=-\frac{5}{2\pi i}\int_{\rho}^{R}
\hat U^{(N)}(t)\,dt + r^{(1)}_N(\rho)  +  r^{(2)}_N(R),
\end{equation}
where the error terms  $r^{(1)}_N(\rho)$ and  $r^{(2)}_N(R)$ satisfy the estimates,
\begin{equation}\label{r1N}
|r^{(1)}_N(\rho)| < C^{(1)}\rho^{5N}, \quad \forall N\geq1,
\end{equation}
and
\begin{equation}\label{r2N}
|r^{(2)}_N(R)| < C_N R^{-5}, \quad \forall R >  R_0,
\end{equation}
respectively. We remind the reader that the constants $C_N$ depends on $N$ only while the constant
$C^{(1)}$  is a numerical constant independent of $N$. In derivation of (\ref{new13})        
 we also took into account symmetries (\ref{u0_gen}) which allowed
us to replace  the sum of integrals from (\ref{residue0_integral_app}) by a single integral.

From (\ref{new6})  it follows that the integrand in (\ref{new13}) satisfies the
asymptotic equation,
 \begin{equation} 
\label{UN_estim1}
\hat U^{(N)}(t)\bigr|_{\arg t=0}=
e^{i\frac{\pi}{2}(k+1)}
\frac{3^{1/4}}{2\sqrt\pi}\,t^{5N-9/4}
e^{-\frac{8\sqrt3}{5}t^{5/2}}
\Bigl(1+r(t)\bigr)\Bigr),
\end{equation}
with
$$
|r(t)| < Ct^{- 5/2}, \quad \forall t : t > R_0.
$$
In particular, this means that the integral of $\hat U^{(N)}(t)$ along the half line  $[\rho, \infty)$ converges. The coefficient $a_N$ and the term $r^{(1)}(\rho)$ in (\ref{new13}) do not depend{\footnote{
Indeed, we have that
$$
r^{(1)}(\rho) = \frac{1}{2\pi i}\oint_{|t|=\rho}\Bigl(\hat u^{(N)}(t) + P_{5N-1}(t)\Bigr) \,dt.
$$}} on
$R$ while the term $r^{(2)}$, in view of (\ref{r2N}), vanishes as $R \to \infty$. Therefore,
sending $R \to \infty$ in (\ref{new13}) (and keeping $N$ fixed) we arrive at the equation,
\begin{equation}\label{aN00}
a_N=-\frac{5}{2\pi i}\int_{\rho}^{\infty}
\hat U^{(N)}(t)\,dt + r^{(1)}_N(\rho),
\end{equation}
which in turn can be transformed as follows.
\begin{multline}\label{Gamma_extraction}
a_N=-\frac{5}{2\pi i}\int_{\rho}^{\infty}
\hat U^{(N)}(t)\,dt
+ r^{(1)}_N(\rho)   =
\\
-\frac{5}{2\pi i}\int_{\rho}^{\infty}
e^{i\frac{\pi}{2}}
\frac{3^{1/4}}{2\sqrt\pi}\,t^{5N-9/4}
e^{-\frac{8\sqrt3}{5}t^{5/2}}\,dt
\\
-\frac{5}{2\pi i}\int_{\rho}^{\infty}
\Bigl(
\hat U^{(N)}(t)
-e^{i\frac{\pi}{2}}
\frac{3^{1/4}}{2\sqrt\pi}\,t^{5N-9/4}
e^{-\frac{8\sqrt3}{5}t^{5/2}}
\Bigr)\,dt
+ r^{(1)}_N(\rho)   =
\\
=-\frac{5\cdot3^{1/4}}{4\pi^{3/2}}
\int_{0}^{\infty}
t^{5N-9/4}
e^{-\frac{8\sqrt3}{5}t^{5/2}}\,dt
\\
+\frac{5\cdot3^{1/4}}{4\pi^{3/2}}
\int_{0}^{\rho}
t^{5N-9/4}
e^{-\frac{8\sqrt3}{5}t^{5/2}}\,dt
-\frac{5}{2\pi i}\int_{\rho}^{\infty}
\Bigl(
\hat U^{(N)}(t)
-e^{i\frac{\pi}{2}}
\frac{3^{1/4}}{2\sqrt\pi}\,t^{5N-9/4}
e^{-\frac{8\sqrt3}{5}t^{5/2}}
\Bigr)\,dt
\\
+ r^{(1)}_N(\rho)   
\equiv I_1+I_2+I_3+
+ r^{(1)}_N(\rho).
\end{multline}

The first integral, $I_1$, in the last  equation  reduces to the Gamma-function integral,
\begin{multline}\label{I1}
I_1=-\frac{5\cdot3^{1/4}}{4\pi^{3/2}}
\int_{0}^{\infty}
t^{5N-9/4}
e^{-\frac{8\sqrt3}{5}t^{5/2}}\,dt=
-\frac{3^{1/4}}{2\pi^{3/2}}
\bigl(
\frac{8\sqrt3}{5}
\bigr)^{-2N+1/2}
\int_{0}^{\infty}
s^{2N-3/2}
e^{-s}\,ds=
\\
=-\frac{3^{1/4}}{2\pi^{3/2}}
\bigl(
\frac{8\sqrt3}{5}
\bigr)^{-2N+1/2}
\Gamma(2N-\tfrac{1}{2}).
\end{multline}
For the second integral, $I_2$, we have
\begin{multline}\label{I2}
|I_2|=\frac{5\cdot3^{1/4}}{4\pi^{3/2}}
\biggl|
\int_{0}^{\rho}
t^{5N-9/4}
e^{-\frac{8\sqrt3}{5}t^{5/2}}\,dt
\biggr|
\leq
\frac{5\cdot3^{1/4}}{4\pi^{3/2}}
\int_{0}^{\rho}
t^{5N-9/4}\,dt=
\\
=\frac{5\cdot3^{1/4}}{4\pi^{3/2}}
\int_{0}^{\rho}
t^{5N-9/4}\,dt=
\frac{3^{1/4}}{4\pi^{3/2}(N-1/4)}
\rho^{5N-5/4} < C^{(2)}\frac{\rho^{5N}}{N}.
\end{multline}
For the third integral, $I_3$, we use the asymptotics (\ref{UN_estim1}),
\begin{equation}\label{I4}
|I_3| < C\frac{5\cdot3^{1/4}}{4\pi \sqrt{\pi}}\int_{\rho}^{\infty}
t^{5N-19/4}
e^{-\frac{8\sqrt3}{5}t^{5/2}}\,dt
< C^{(3)}\Gamma(2N-\tfrac{3}{2}) + C^{(4)}\frac{\rho^{5N}}{N},
\end{equation}
where the constants $C^{(j)}$, $j =2, 3, 4$ are positive constants which do not depend on $N$.
Equation (\ref{I1}) and estimates (\ref{I2}), (\ref{I4}), and (\ref{r1N}) mean that
\begin{equation}\label{aN0}
a_N = -\frac{3^{1/4}}{2\pi^{3/2}}
\bigl(
\frac{8\sqrt3}{5}
\bigr)^{-2N+1/2}
\Gamma(2N-\tfrac{1}{2})
\bigl(1+{\mathcal O}(N^{-1})\bigr) + {\mathcal O}(\rho^{5N}),
\end{equation}
as $N \to \infty$. This, in turn, implies the following expression
for the leading term of the large $N$ asymptotics of the coefficient $a_N$.
\begin{equation}\label{residue0_integral}
a_N = =-\frac{3^{1/4}}{2\pi^{3/2}}
\bigl(
\frac{8\sqrt3}{5}
\bigr)^{-2N+1/2}
\Gamma(2N-\tfrac{1}{2})
\bigl(1+{\mathcal O}(N^{-1})\bigr), \quad N \to \infty.
\end{equation}

Equation \eqref{residue0_integral} provides us with the leading term of the
large $n$ asymptotics of $a_n \equiv u_{n,0}$. In order to reproduce the
whole series \eqref{aas} we note that the Riemann-Hilbert analysis
of the tritronqu\'ee solutions yields in fact the full asymptotic series
expansion for the differences between the tritronqu\'ee solutions.
This means the following extension of the estimates \eqref{quasi_linear_Stokes}:
\begin{equation}\label{extention1}
u_{k+1}(z)-u_k(z)\sim e^{i\frac{\pi}{2}(k+1)}
\frac{3^{1/4}}{2\sqrt\pi}\,z^{-1/8}
e^{(-1)^{k+1}\frac{8\sqrt3}{5}z^{5/4}}
\Bigl(1+\sum_{n=1}^{\infty}c_nz^{-5n/4}\Bigr),
\end{equation}
$$
z\to\infty,\quad
\arg z\in\bigl(-\tfrac{2\pi}{5}+\tfrac{4\pi}{5}k,
\tfrac{2\pi}{5}+\tfrac{4\pi}{5}k\bigr),
$$
with some coefficients $c_n$. In its turn, asymptotics \eqref{extention1} implies that 
\begin{equation}\label{extension2}
\hat U^{(N)}(t)\bigr|_{\arg t=\tfrac{2\pi}{5}k}=
e^{i\frac{\pi}{2}(k+1)}
\frac{3^{1/4}}{2\sqrt\pi}\,t^{5N-9/4}
e^{(-1)^{k+1}\frac{8\sqrt3}{5}t^{5/2}}
\Bigl(1+\sum_{n=1}^{m}c_nt^{-5n/2} +{\mathcal O}(t^{-5(m+1)/2})\Bigr).
\end{equation}
Therefore, the integral of $\hat U^{(N)}(t)$ along the half-line $(\rho, \infty)$ can
be estimated {\footnote{ We omit  the routine technical details which are similar to
the ones carefully presented in the 
derivation of (\ref{residue0_integral} )}}  as follows:
\be
\ba
\label{extension3}
& \int_{\rho}^{\infty} \hat U^{(N)}(t)\,dt = -S_1\left(\sum_{n=0}^{m}c_n\int_{\rho}^{\infty}
e^{-At^{5/2}}t^{5N-5n/2 -9/4}\,dt + {\mathcal O}\Bigl(\int_{\rho}^{\infty}
e^{-At^{5/2}}t^{5N-5(m+1)/2 -9/4}\,dt\Bigr)\right)\\
&= -\frac{2}{5}S_1\left(\sum_{n=0}^{m}c_nA^{-2N+n+1/2}\int_{A\rho^{5/2}}^{\infty}
e^{-t}t^{2N-n -3/2}\,dt + {\mathcal O}\Bigl(A^{-2N+m+3/2}\int_{A\rho^{5/2}}^{\infty}
e^{-t}t^{2N-m-1 - 3/2}\,dt\Bigr)\right)\\ 
&= -\frac{2}{5}S_1\left(\sum_{n=0}^{m}c_nA^{-2N+n+1/2}\Gamma(2N-n-\tfrac{1}{2})
+ {\mathcal O}\Bigl(A^{-2N+m+3/2}\Gamma(2N - m-\tfrac{3}{2})\Bigr)\right)
+ {\mathcal O} (\rho^{5N -5/4})\\
&= -\frac{2}{5}A^{-2N+1/2}\Gamma(2N-\tfrac{1}{2})S_1\left(1 + \sum_{n=1}^{m}
\frac{c_nA^{n}}{\prod_{k=1}^{n}\left(2N-1/2 -k\right)}
+ {\mathcal O}\left(\frac{A^{m+1}}{\prod_{k=1}^{m+1}\left(2N-1/2 -k\right)}\right)\right)
+ {\mathcal O} (\rho^{5N -5/4})\\
&= -\frac{2}{5}A^{-2N+1/2}\Gamma(2N-\tfrac{1}{2})S_1\left\{1 + \sum_{n=1}^{m}
\frac{c_nA^{n}}{\prod_{k=1}^{n}\left(2N-1/2 -k\right)}
+ {\mathcal O}\left(N^{-m-1}\right)\right\}.
\ea
\ee
Here, $m$ is arbitrary but fixed and all the constants in the error terms depend on $m$ and $\rho$
(which is also fixed) only. We also have used our standard  notations,
$$
A = \frac{8\sqrt{3}}{5}, \quad S_1 = -\ri {3^{1\over 4} \over 2{\sqrt{\pi}}}, \quad c_0=1.
$$
Substituting \eqref{extension3} into \eqref{aN00} we arrive at
the final estimate for $a_N$,
\begin{equation}\label{aNfinal}
a_N = A^{-2N+1/2}\Gamma(2N-\tfrac{1}{2})\frac{S_1}{\pi \ri}\left\{1 + \sum_{n=1}^{m}
\frac{c_nA^{n}}{\prod_{k=1}^{n}\left(2N-1/2 -k\right)}
+ {\mathcal O}\left(N^{-m-1}\right)\right\}.
\end{equation}
To complete the proof of the $0$ - instanton asymptotics \eqref{aas} we need to identify 
the coefficients $c_n$ as the $1$ -instanton coefficients $u_{n,1}$. To this end, it is enough to notice
that the difference, $v(z) \equiv u_{k+1}(z) - u_{k}(z)$, of the tritronqu\'ee solutions satisfies 
the linear differential equation,
\begin{equation}\label{extension5}
 v'' = 6v(u_{k+1} + u_{k}).
 \end{equation}
 Since both $u_{k}$ and $u_{k+1}$ have the same power series expansion in the
 relevant sectors, we conclude that the formal power series solution of equation
 \eqref{extension5} must coincide with the formal power-series solution of
 the $1$ -instanton equation \eqref{one000}. Hence the desired equation,
 \begin{equation}\label{extension6}
 c_n = u_{n,1}, \quad \forall n,
 \end{equation}
which completes the proof of \eqref{aas}.

\sectiono{Asymptotics of the $1$-instanton of Painlev\'e I}
\label{sec.1inst}

\subsection{A 5-tuple of differential equations}
Let us consider the first instanton term $Cu_1$ in the formal series 
\eqref{uzC}. This term satisfies the linear homogeneous second order ODE
\begin{equation}
\label{1instanton_ODE}
u_1''=12u_0u_1,
\end{equation}
where $u_0$ is the 0-instanton term in the trans-series expansion 
\eqref{uzC}. This suggests the problem of studying  the linear equation
\begin{equation}
\label{gamma_ODE}
v''=\gamma u v,
\end{equation}
where $u$ is one of the tritronqu\'ee solutions of the first Painlev\'e equation
\eqref{PI} and $\gamma$ is an arbitrary  positive constant. 
Indeed, we have five {\it different} Schr\"odinger
equations with meromorphic potentials $\gamma u_j$ 
having the {\it same}  formal 
power expansion \eqref{u0_formal} in the relevant sectors, see \eqref{u0_gen}:
\begin{equation}
\label{potential_as}
v_j''=\gamma u_j v_j,\quad
\gamma u_j\sim
\gamma u_f(z)=\gamma z^{1/2}\sum_{n=0}^{\infty}a_n z^{-5n/2}.
\end{equation}
Equation \eqref{gamma_ODE} has two linearly independent 
formal solutions,
\begin{multline}
\label{v_formal}
v_f^{\pm}=e^{\theta^{\pm}}\Sigma^{\pm},\quad
\theta^{\pm}= B^{\pm} z^{5/4}-\tfrac{1}{8}\ln z,\quad
\Sigma^{\pm}=\sum_{n=0}^{\infty}b_n^{\pm}z^{-5n/4},
\\
B^{\pm}=\pm B=\pm\tfrac{4}{5}\sqrt{\gamma},\quad
b_0^{\pm}=1,\quad
b_n^{\pm}=\frac{1}{2n}
\Bigl\{
\frac{b_{n-1}^{\pm}}{B^{\pm}}\bigl(n-\tfrac{1}{10}\bigr)
\bigl(n-\tfrac{9}{10}\bigr)
-B^{\pm}\sum_{m=1}^{[\frac{n+1}{2}]}a_mb_{n+1-2m}^{\pm}
\Bigr\},\quad
n=1,2,\dots
\end{multline}
To simplify our notation, in this section, we use the symbol $b_n^-$ instead 
of $u_{n,1}$. 

Since we have five equations with meromorphic potentials, we have five 
pairs of solutions with the asymptotic expansions \eqref{v_formal} as $z\to\infty$.
Each pair has this asymptotic expansion in the sector where the
corresponding potential has the power-like asymptotics \eqref{u0_formal}.
The sectors are depicted  in  \figref{tritronquee}.
Denote the solutions corresponding to the potential $u_j(z)$ by the 
symbol $v_j^{\pm}(z)$. We have,
\begin{equation}
\label{vj_gen}
v_j^{\pm}(z)
\sim
v_f^{\pm}(z), \quad
z\to\infty,\quad
\arg z\in\bigl(-\tfrac{6\pi}{5}+\tfrac{4\pi}{5}j,\tfrac{2\pi}{5}
+\tfrac{4\pi}{5}j\bigr).
\end{equation}
We also observe that,
\begin{equation}
\label{vj_gen2}
v_j^{\pm}(z)=e^{-i\frac{\pi}{10}j}v_0^{\sigma_j(\pm)}(e^{-i\frac{4\pi}{5}j}z),
\quad
\sigma_j(\pm)=
\begin{cases}
\pm,\quad
\mbox{$j$ is even},\\
\mp,\quad
\mbox{$j$ is odd},
\end{cases}
\end{equation}
and
\begin{equation}\label{vj_gen3}
v_{j+20}^{\pm}(z) = v_j^{\pm}(z).
\end{equation}

Similar to convention \eqref{module} concerning the subscribe $k$ in $u_k(x)$,
we shall assume that
\begin{equation}\label{module2}
j \in {\Bbb Z},\quad \mbox{mod}\quad 20,
\end{equation}
in the notation $ v_j^{\pm}(z)$, unless $j$ is particularly specified, as in \eqref{hat_V_def}
below.
 
\subsection{Gluing the solutions together}

We observe that within the common sectors of the same asymptotic behavior,
the solutions $v_{j}^{\pm}(z)$ and $v_{j-1}^{\pm}(z)$ have identical 
asymptotics
$v_f^{\pm}(z)$ in all orders. Thus it is natural to ask: what is the 
difference
between these two solutions of {\em different} equations? Let us introduce 
the 
differences
\begin{equation}
\label{w_def}
w_j^{\pm}(z)=v_{j+1}^{\pm}(z)-v_{j}^{\pm}(z).
\end{equation}
Because of \eqref{vj_gen2}, these differences are related to each other by the 
following rotational symmetry,
\begin{equation}
\label{wj_w0}
w_j^{\pm}(z)=e^{-i\frac{\pi}{10}j}w_0^{\sigma_j(\pm)}(e^{-i\frac{4\pi}{5}j}z).
\end{equation}
In addition, using \eqref{potential_as}, $w_j^{\pm}(z)$ satisfy 
non-homogeneous 
linear ODE,
\begin{equation}
\label{w_eq}
(w_{j-1}^{\pm})_{zz}=
\gamma u_{j}w_{j-1}^{\pm}
+\gamma\bigl(u_{j}-u_{j-1}\bigr)v_{j-1}^{\pm}.
\end{equation}
The homogeneous part of this equation is the same as in \eqref{potential_as}
In the relevant sector, the non-homogeneity contains an exponentially 
small 
factor given in \eqref{quasi_linear_Stokes}:
\begin{multline}
\label{induced_linear_Stokes}
z\to\infty,\quad
\arg z\in\bigl(-\tfrac{6\pi}{5}+\tfrac{4\pi}{5}j,
-\tfrac{2\pi}{5}+\tfrac{4\pi}{5}j\bigr)\colon
\\
u_{j}(z)-u_{j-1}(z)=e^{i\frac{\pi}{2}j}
\frac{3^{1/4}}{2\sqrt\pi}\,
e^{\theta_j}
\Bigl(1+{\mathcal O}\bigl(z^{-5/4}\bigr)\Bigr),
\\
\theta_j=\alpha_jz^{5/4}-\frac{1}{8}\ln z,\quad
\alpha_j=(-1)^jA,\quad
A=\frac{8\sqrt3}{5}.
\end{multline}

The general solution to \eqref{w_eq} is given by the integral formula,
\begin{equation}
\label{w_general}
w_{j-1}^{\pm}(z)=c^+v_j^+(z)+c^-v_j^-
-\frac{1}{2\sqrt{\gamma}}
\int_{z_0}^z\bigl(v_j^+(x)v_j^-(z)-v_j^+(z)v_j^-(x)\bigr)
\gamma\bigl(u_{j}(x)-u_{j-1}(x)\bigr)v_{j-1}^{\pm}(x)\,dx,
\end{equation}
where $c^+$, $c^-$ and $z_0$ are arbitrary constants. However, our 
solution \eqref{w_def} is {\em not} general as being a difference between
two functions with the same asymptotic expansion in certain sector. 
Furthermore, in the sector 
$\arg z\in\bigl(-\tfrac{6\pi}{5}+\tfrac{4\pi}{5}j,
-\tfrac{2\pi}{5}+\tfrac{4\pi}{5}j\bigr)$, see \eqref{induced_linear_Stokes},
for $j$ odd, the solution $v_j^+(z)$ is dominant and $v_j^-(z)$ is
recessive, while for $j$ even, $v_j^+(z)$ is recessive and $v_j^-(z)$ is
dominant. Thus $w_{j-1}^{\pm}(z)$ admits the following representations:
\begin{subequations}
\label{w_representations}
\begin{multline}
\label{w_rep_j_odd}
\mbox{$j$ odd: $v_j^+(z)$ dominant, $v_j^-(z)$ recessive as 
$\arg z\in\bigl(-\tfrac{6\pi}{5}+\tfrac{4\pi}{5}j,
-\tfrac{2\pi}{5}+\tfrac{4\pi}{5}j\bigr)$}\colon
\\
\shoveleft
w_{j-1}^+(z)=c_0v_j^-(z)
-\frac{1}{2}\sqrt{\gamma}
v_j^-(z)
\int_{e^{i\tfrac{4\pi}{5}(j-1)}z_0}^z
\bigl(u_{j}(x)-u_{j-1}(x)\bigr)
v_j^+(x)v_{j-1}^+(x)\,dx
\\
\hfill
+\frac{1}{2}\sqrt{\gamma}
v_j^+(z)
\int_{e^{i\tfrac{4\pi}{5}(j-1)}\infty}^z
\bigl(u_{j}(x)-u_{j-1}(x)\bigr)
v_j^-(x)v_{j-1}^+(x)\,dx,
\\
\shoveleft
w_{j-1}^-(z)=
-\frac{1}{2}\sqrt{\gamma}
\int_{e^{i\tfrac{4\pi}{5}(j-1)}\infty}^z
\bigl(v_j^+(x)v_j^-(z)-v_j^+(z)v_j^-(x)\bigr)
\bigl(u_{j}(x)-u_{j-1}(x)\bigr)v_{j-1}^-(x)\,dx,
\hfill
\end{multline}
\begin{multline}
\label{w_rep_j_even}
\mbox{$j$ even: $v_j^+(z)$ recessive, $v_j^-(z)$ dominant as 
$\arg z\in\bigl(-\tfrac{6\pi}{5}+\tfrac{4\pi}{5}j,
-\tfrac{2\pi}{5}+\tfrac{4\pi}{5}j\bigr)$}\colon
\\
\shoveleft
w_{j-1}^+(z)=
-\frac{1}{2}\sqrt{\gamma}
\int_{e^{i\tfrac{4\pi}{5}(j-1)}\infty}^z
\bigl(v_j^+(x)v_j^-(z)-v_j^+(z)v_j^-(x)\bigr)
\bigl(u_{j}(x)-u_{j-1}(x)\bigr)v_{j-1}^+(x)\,dx,
\\
\shoveleft
w_{j-1}^-(z)=c_0v_j^+(z)
-\frac{1}{2}\sqrt{\gamma}
v_j^-(z)
\int_{e^{i\tfrac{4\pi}{5}(j-1)}\infty}^z
\bigl(u_{j}(x)-u_{j-1}(x)\bigr)
v_j^+(x)v_{j-1}^-(x)\,dx
\\
+\frac{1}{2}\sqrt{\gamma}
v_j^+(z)
\int_{e^{i\tfrac{4\pi}{5}(j-1)}z_0}^z
\bigl(u_{j}(x)-u_{j-1}(x)\bigr)
v_j^-(x)v_{j-1}^-(x)\,dx,
\end{multline}
\end{subequations}
where $e^{i\tfrac{4\pi}{5}(j-1)}z_0$ is a finite point within the indicated 
sector and $c_0$ is a Stokes constant which both can not be determined immediately.

When $z\to\infty$ as $\arg z\in\bigl(-\tfrac{6\pi}{5}+\tfrac{4\pi}{5}j,
-\tfrac{2\pi}{5}+\tfrac{4\pi}{5}j\bigr)$, the leading contribution to
the asymptotic behavior of $w_{j-1}^{\pm}(z)$ is given by $z$-end point of the 
integral term and possibly by the non-integral term. By standard arguments,
\begin{subequations}
\label{w_as}
\begin{multline}
\label{w_as_j_odd}
\mbox{$j$ odd: $v_j^+(z)$ dominant, $v_j^-(z)$ recessive as 
$\arg z\in\bigl(-\tfrac{6\pi}{5}+\tfrac{4\pi}{5}j,
-\tfrac{2\pi}{5}+\tfrac{4\pi}{5}j\bigr)$}\colon
\\
\\
\shoveleft
w_{j-1}^+(z)=\begin{cases}
\displaystyle{-e^{i\frac{\pi}{2}j}
\frac{3^{1/4}}{5\sqrt\pi}
\frac{2B\sqrt{\gamma}}{A(-A+2B)}
e^{(-A+B)z^{5/4}}z^{-3/4}
\bigl(1+{\mathcal O}(z^{-5/4})\bigr),\quad
B>A/2},\\
\\
\displaystyle{-e^{i\frac{\pi}{2}j}
\frac{2\cdot3^{1/4}}{5\sqrt\pi}\,
\sqrt{\gamma}
e^{-B z^{5/4}}z^{1/2}
\bigl(1+{\mathcal O}(z^{-5/8})\bigr),\quad
B=A/2},\\
\\
\displaystyle{c_1e^{-B z^{5/4}}z^{-1/8}
\bigl(1+{\mathcal O}(z^{-5/4})\bigr),\quad
0<B<A/2},
\end{cases}
\\
\\
\shoveleft
w_{j-1}^-(z)=
e^{i\frac{\pi}{2}j}
\frac{3^{1/4}}{5\sqrt\pi}
\frac{2B\sqrt{\gamma}}{A(A+2B)}
e^{-(A+B) z^{5/4}}z^{-3/4}
\bigl(1+{\mathcal O}(z^{-5/4})\bigr),
\hfill
\end{multline}
\begin{multline}
\label{w_as_j_even}
\mbox{$j$ even: $v_j^+(z)$ recessive, $v_j^-(z)$ dominant as 
$\arg z\in\bigl(-\tfrac{6\pi}{5}+\tfrac{4\pi}{5}j,
-\tfrac{2\pi}{5}+\tfrac{4\pi}{5}j\bigr)$}\colon
\\
\\
\shoveleft
w_{j-1}^+(z)=
e^{i\frac{\pi}{2}j}
\frac{3^{1/4}}{5\sqrt\pi}
\frac{2B\sqrt{\gamma}}{A(A+2B)}
e^{(A+B)z^{5/4}}z^{-3/4}
\bigl(1+{\mathcal O}(z^{-5/4})\bigr)
\\
\\
\shoveleft
w_{j-1}^-(z)=
\begin{cases}
\displaystyle{e^{i\frac{\pi}{2}j}
\frac{3^{1/4}}{5\sqrt\pi}
\frac{2B\sqrt{\gamma}}{A(A-2B)}
e^{(A-B)z^{5/4}}z^{-3/4}
\bigl(1+{\mathcal O}(z^{-5/4})\bigr),\quad
B>A/2},
\\
\\
\displaystyle{e^{i\frac{\pi}{2}j}
\frac{2\cdot3^{1/4}}{5\sqrt\pi}
\sqrt{\gamma}
e^{B z^{5/4}}z^{1/2}
\bigl(1+{\mathcal O}(z^{-5/8})\bigr),\quad
B=A/2},\\
\\
\displaystyle{c_1e^{B z^{5/4}}z^{-1/8}
\bigl(1+{\mathcal O}(z^{-5/4})\bigr),\quad
0<B<A/2},\\
\end{cases}
\hfill
\end{multline}
\end{subequations}
where $c_1$ is an unknown Stokes constant.

Taking into account the presence of the factor $z^{-8}$ in the asymptotics
of the functions $ v_j^{\pm}(\xi)$,  we define  auxiliary functions $\hat v_j^{\pm}(\xi)$
by the relations (cf. \eqref{hat_uk_def}),
\begin{equation}
\label{hat_vj_def}
z=\xi^8,\quad
\hat v_j^{\pm}(\xi)=v_j^{\pm}(z)e^{-\theta^{\pm}(z)},
\end{equation}
so that (see (\ref{vj_gen}))
\begin{equation}
\label{hat_vj_as}
\hat v_j^{\pm}(\xi)\sim
\Sigma^{\pm}(\xi^8)=\sum_{n=0}^{\infty}b_n^{\pm}\xi^{-10n},\quad
\xi\to\infty,\quad
\arg\xi\in\bigl(-\tfrac{3\pi}{20}+\tfrac{\pi}{10}j,\tfrac{\pi}{20}
+\tfrac{\pi}{10}j\bigr).
\end{equation}
Using \eqref{vj_gen2}, we observe that
\begin{equation}
\label{hat_vj_symm}
\hat v_j^{\pm}(\xi)=e^{-i\frac{\pi}{10}j}
\hat v_0^{\sigma_j(\pm)}\bigl(e^{-i\frac{\pi}{10}j}\xi\bigr),\quad
\hat v_{j+20}^{\pm}(\xi)=\hat v_j^{\pm}(\xi).
\end{equation}
Let us introduce the sectorially meromorphic function $\hat V^{(N)}(\xi)$ 
discontinuous
across the rays $\arg\xi=\frac{\pi}{10}j$, $j=-9,-8,\dots,9,10$:
\begin{equation}
\label{hat_V_def}
\hat V^{(N)}(\xi)=\hat v_j^-(\xi)\xi^{10N-1}-P_{10N-1}^-(\xi),\quad
\arg\xi\in\bigl(-\tfrac{\pi}{10}+\tfrac{\pi}{10}j,\tfrac{\pi}{10}j\bigr),\quad
j=-9,-8,\dots,9,10,
\end{equation}
where polynomial $P_{10N-1}^-(\xi)$ of degree $10N-1$ is defined by
\begin{equation}
\label{P10N-1_def}
P_{10N-1}^-(\xi):=\sum_{n=0}^{N-1}b_n^-\xi^{10(N-n)-1}.
\end{equation}
According to \eqref{hat_vj_as}, this function has the uniform asymptotics
at infinity,
\begin{equation}
\label{VN_as}
V^{(N)}(\xi)=b_n^-\xi^{-1}
+{\mathcal O}_{N}(\xi^{-11}),\quad
\xi\to\infty,
\end{equation}
and, due to \eqref{w_def} and \eqref{hat_vj_def}, has the following jumps 
across the rays $\arg\xi=\frac{\pi}{10}(j-1)$, $j=-9,-8,\dots,9,10$, oriented 
towards infinity:
\begin{multline}
\label{hat_V_jumps}
\arg\xi=\tfrac{\pi}{10}(j-1)\colon\quad
\hat V^{(N)}_+(\xi)-\hat V^{(N)}_-(\xi)=
\bigl(\hat v_j^-(\xi)-\hat v_{j-1}^-(\xi)\bigr)\xi^{10N-1}=
\\
=\bigl(v_j^-(z)-v_{j-1}^-(z)\bigr)e^{-\theta^-(z)}\xi^{10N-1}=
w_{j-1}^-(z)e^{-\theta^-(z)}\xi^{10N-1},\quad
j=-8,-7,\dots,11.
\end{multline}
(We note that $v^{-}_{11}(z) = v^{-}_{-9}(z)$,
in virtue of \eqref{vj_gen3}). The jumps for even and odd $j$-s are different.
Namely, as $\xi\to\infty$, using \eqref{w_as}, we find
\begin{subequations}
\label{hat_V_jumps_as}
\begin{multline}
\label{hat_V_jumps_as_j_odd}
\xi\to\infty,\quad
\arg\xi=\tfrac{\pi}{10}(j-1),\quad
\mbox{$j$ odd}\colon
\\
\shoveleft
{
\hat V^{(N)}_+(\xi)-\hat V^{(N)}_-(\xi)=
e^{i\frac{\pi}{2}j}
\frac{3^{1/4}}{5\sqrt\pi}
\frac{2B\sqrt{\gamma}}{A(A+2B)}
e^{-A \xi^{10}}\xi^{10N-6}
\bigl(1+{\mathcal O}(\xi^{-10})\bigr),
}
\hfill
\end{multline}
\begin{multline}
\label{hat_V_jumps_as_j_even}
\xi\to\infty,\quad
\arg\xi=\tfrac{\pi}{10}(j-1),\quad
\mbox{$j$ even}\colon
\\
\shoveleft
{
\hat V^{(N)}_+(\xi)-\hat V^{(N)}_-(\xi)=
\begin{cases}\displaystyle{
e^{i\frac{\pi}{2}j}
\frac{3^{1/4}}{5\sqrt\pi}
\frac{2B\sqrt{\gamma}}{A(A-2B)}
e^{A\xi^{10}}\xi^{10N-6}
\bigl(1+{\mathcal O}(\xi^{-10})\bigr),\quad
B>A/2},\\
\\
\displaystyle{e^{i\frac{\pi}{2}j}
\frac{2\cdot3^{1/4}}{5\sqrt\pi}
\sqrt{\gamma}
e^{A\xi^{10}}\xi^{10N+4}
\bigl(1+{\mathcal O}(\xi^{-5})\bigr),\quad
B=A/2},\\
\\
\displaystyle{c_1e^{2B\xi^{10}}\xi^{10N-1}
\bigl(1+{\mathcal O}(\xi^{-10})\bigr),\quad
0<B<A/2}.
\end{cases}
}
\hfill
\end{multline}
\end{subequations}

\begin{remark}\label{routine}
In this section, as well as in the rest of the paper, we skip the detail explanation
of the standard meaning of the symbols ${\mathcal O}(...)$ and  ${\mathcal O}_{N}(...)$.
This meaning has been spelled out in Section \ref{as.0inst} (see e.g. (\ref{new3}),
(\ref{new4}), (\ref{new6}), and (\ref{new10})). We also omit, when performing 
asymptotic calculations, the routine technical details. They are similar to the ones
which were carefully presented in Section \ref{as.0inst}.
\end{remark}  

\subsection{An integral formula for the $1$-instanton coefficients and their asymptotic
expansion}
\lbl{sub.integral}

Due to \eqref{VN_as}, the 1-instanton coefficient $b_N^-$, whose $N$-large 
asymptotics 
we are looking for, is the coefficient at the term $1/\xi$  of the asymptotic
expansion of $V^{(N)}(\xi)$ at $\xi = \infty$. Therefore,
\begin{equation}
\label{residue1_infinity}
\frac{1}{2\pi i}\oint_{|\xi|=R\gg1}V^{(N)}(\xi)\,d\xi
=b_n^-+{\mathcal O}_{N}(R^{-10}).
\end{equation}
Similarly to the $0$-instanton case, collapsing this circular path of 
integration 
to a circle of a smaller radius $|\xi|=r$ containing all the singularities 
or branch 
points of $\hat V^{(N)}(\xi)$ (the latter can appear at the poles of 
$u_j(\xi^8)$ only),
we find
\begin{multline}
\label{residue1_itegral}
b_N^-=\frac{1}{2\pi i}\oint_{|\xi|=R\gg1}V^{(N)}(\xi)\,d\xi
+{\mathcal O}_{N}(R^{-10})=
\\
=-\frac{1}{2\pi i}\sum_{j=-9}^{10}
\int_{re^{i\frac{\pi}{10}(j-1)}}^{Re^{i\frac{\pi}{10}(j-1)}}
w_{j-1}^-(z)e^{-\theta^-(z)}\xi^{10N-1}\,d\xi
+\frac{1}{2\pi i}\oint_{|\xi|=r}V^{(N)}(\xi)\,d\xi
+{\mathcal O}_{N}(R^{-10})=
\\
=-\frac{1}{2\pi i}\sum_{j=-9}^{10}
\int_{re^{i\frac{\pi}{10}(j-1)}}^{Re^{i\frac{\pi}{10}(j-1)}}
w_{j-1}^-(z)e^{-\theta^-(z)}\xi^{10N-1}\,d\xi
+\frac{1}{2\pi i}\oint_{|\xi|=r}\Bigl(V^{(N)}(\xi)+P_{10N-1}^{-}(\xi)\Bigr)\,d\xi
+{\mathcal O}_{N}(R^{-10})=
\\
=-\frac{1}{2\pi i}\sum_{j=-9}^{10}
\int_{re^{i\frac{\pi}{10}(j-1)}}^{Re^{i\frac{\pi}{10}(j-1)}}
e^{-i\frac{\pi}{10}(j-1)}w_0^{\sigma_{j-1}(-)}(e^{-i\frac{4\pi}{5}(j-1)}z)
e^{-\theta^-(z)}\xi^{10N-1}\,d\xi
+{\mathcal O}(r^{10N})
+{\mathcal O}_{N}(R^{-10})=
\\
=-\frac{1}{2\pi i}\sum_{k=-5}^{4}
\int_{re^{i\frac{\pi}{5}k}}^{Re^{i\frac{\pi}{5}k}}
e^{-i\frac{\pi}{5}k}w_0^-\bigl((e^{-i\frac{\pi}{5}k}\xi)^8\bigr)
e^{B\xi^{10}}\xi^{10N}\,d\xi
\\
-\frac{1}{2\pi i}\sum_{k=-5}^{4}
\int_{re^{i\frac{\pi}{10}(2k+1)}}^{Re^{i\frac{\pi}{10}(2k+1)}}
e^{-i\frac{\pi}{10}(2k+1)}w_0^+\bigl((e^{-i\frac{\pi}{10}(2k+1)}\xi)^8\bigr)
e^{B\xi^{10}}\xi^{10N}\,d\xi
+{\mathcal O}(r^{10N})
+{\mathcal O}_{N}(R^{-10})=
\\
=-\frac{5}{\pi i}
\int_{r}^{R}
w_0^-\bigl(\xi^8\bigr)
e^{B\xi^{10}}
\xi^{10N}\,d\xi
-\frac{5}{\pi i}
(-1)^N
\int_{r}^{R}
w_0^+\bigl(\xi^8\bigr)
e^{-B\xi^{10}}
\xi^{10N}\,d\xi
+{\mathcal O}(r^{10N})
+{\mathcal O}_{N}(R^{-10}).
\end{multline}
Here, we  used that
$$
|V^{(N)}(\xi)+P_{10N-1}^{-}(\xi)|_{|\xi| = r} < C r^{10N-1},
$$
where the constant $C$ does not depend on $N$.
It is not difficult to find $N$-large asymptotics of the coefficients $b_N^-$.
Similar to the $0$ - instanton  case considered in Section \ref{as.0inst},
we first send $R\to \infty$ in the 
last equation of  (\ref{residue1_itegral}) and then, using \eqref{w_as_j_odd}, analyze 
the large $N$ behavior of the resulting integrals:

\begin{equation}
\label{instanton1_as}
b_N^-=
\begin{cases}
\displaystyle{\frac{4\cdot3^{1/4}}{25\pi^{3/2}}
\frac{A(1+(-1)^N)-2B(1-(-1)^N)}{4B^2-A^2}
\gamma
A^{-N-\frac{1}{2}}
\Gamma(N-\tfrac{1}{2})
\bigl(1+{\mathcal O}(N^{-1})\bigr),\quad
\gamma>3},\\
\\
\displaystyle{(-1)^N
\frac{3^{1/4}}{5\pi^{3/2}}\,
\sqrt{\gamma}
(2B)^{-N-\frac{1}{2}}
\Gamma(N+\frac{1}{2})
\bigl(1+{\mathcal O}(N^{-1/2})\bigr),\quad
\gamma=3},\\
\\
\displaystyle{i(-1)^N
c_1
\frac{1}{2\pi}
(2B)^{-N}
\Gamma(N)
\bigl(1+{\mathcal O}(N^{-1})\bigr),\quad
0<\gamma<3}.
\end{cases}
\end{equation}
We recall that $B=\frac{4}{5}\sqrt\gamma$ and $A=\frac{8}{5}\sqrt3$.

Equation \eqref{instanton1_as}, in the case $\gamma = 12$, produces the leading term of the 
asymptotics of the $1$-instanton sequence $u_{n,1}$. The proof of the whole series 
\eqref{ueven1}--\eqref{uodd1} is similar to the proof of the  $0$-instanton
asymptotic series \eqref{aas} performed in the end of Section \ref{as.0inst}.

\begin{remark}
\label{no_Stokes_condition}
Note, that for $\gamma>3$, the unknown Stokes constant $c_1$  does not contribute
to the leading term of the large $N$ asymptotics of $b_{N}^{-}$; indeed, it
appears in the exponentially small terms only.
For $\gamma=3$, this contribution appears in the subleading term of order
$N^{-1/2}$. However, for $0<\gamma<3$, the asymptotics of $b_N^-$ is simply
proportional to $c_1$. This is also confirmed by numerical computations.
\end{remark}

\subsection{Tritronqu\'ee solutions and the induced Stokes' phenomenon}
\lbl{inducedStokes}
The  analysis performed in the previous subsections gives rise to the following
observations concerning the Schr\"odinger equation,
\begin{equation}\label{indStokes1}
v{''} = \gamma uv.
\end{equation}

 Let us take two  different tritronqu\'ee Painlev\'e I functions, 
say the functions  $u_0(z)$ and $u_1(z)$,
as the potentials $u(z)$ in equation \eqref{indStokes1}. The corresponding 
solutions, $v_0^{+}(z)$ and $v_1^{+}(z)$, though solutions to {\it different}
equations, would have the {\it same}  asymptotics in all orders in the common sector,
$-2\pi/5 < \arg z < 2\pi/5 $. The same is true for the pair $v_0^{-}(z)$ and $v_1^{-}(z)$.
Hence the possibility arises of evaluating the exponentially small differences, $v_1^{+}(z) - v_1^{+}(z)$ and 
 $v_1^{-}(z) - v_1^{-}(z)$. In other words, the existence of the tritronqu\'ee Painlev\'e transcendents
allows to introduce the notion of the Stokes' phenomenon
for the collection of the solutions to the family of linear equations \eqref{indStokes1} consisting 
of five equations generated by five tritronqu\'ee functions  as potentials $u(z)$. We 
shall call this phenomenon the {\it induced Stokes' phenomenon} or, more lengthy,  the {\it 
induced Stokes' phenomenon  generated 
by the first Painlev\'e  quasi-linear Stokes' phenomenon}. It is worth noticing
that no single equation \eqref{indStokes1} with $u(z) = u_k(z)$ generates 
a meaningful Stokes phenomenon in the set of its solutions. Indeed, the structure
of the solutions to equation \eqref{indStokes1} in the complementary sector,
where the potential has infinitely many poles (see Figure \ref{tritronquee}), is 
extremely difficult to describe. We need five different potentials in  \eqref{indStokes1}
to cover the neighborhood of infinity by the sectors with the regular behavior of both
potentials and the solutions. 

There is a threshold value of the parameter $\gamma$ in \eqref{indStokes1}, namely,
$\gamma = 3$. When $\gamma > 3$, the induced Stokes' phenomenon generated
by the first Painlev\'e quasi-linear Stokes phenomenon is completely described 
by the latter. Indeed, the origin of the coefficient $3^{1/4}/2\sqrt{\pi}$ in formulae 
\eqref{w_as_j_odd} and \eqref{w_as_j_even} is the Stokes' constant
$$
S_1 = -\ri {3^{1\over 4} \over 2{\sqrt{\pi}}},
$$
which is featuring in the quasi-linear Stokes' relations \eqref{quasi_linear_Stokes}
for the tritronqu\'ee solutions $u_k(z)$. When $\gamma < 3$, the quasi-linear 
Stokes' phenomenon for potentials of  equation \eqref{indStokes1} does not control 
the induced Stokes' phenomenon  for its solutions. The intrinsic Stokes' constant -
the constant $c_1$, begins to play a dominant role.

There are at least two values of the parameter $\gamma$ which are less than
$3$ but  for which we believe an explicit description of the induced Stokes'
phenomenon is possible. These values are{\footnote{ Painlev\'e I equation (\ref{PI}) provides
the linear equation (\ref{gamma_ODE}) with the potential $u(z)$. We note that the value of the parameter $\gamma$
in equation  (\ref{gamma_ODE})  can not be changed via the scaling of the variables 
in the Painlev\'e  equation without violation of the chosen form of the latter. The values (\ref{exceptions})
corresponds to form  (\ref{PI}) of the equation PI chosen in this paper.}},
\begin{equation}\label{exceptions}
\gamma = 2, \quad \mbox{and}\quad \gamma = \frac{3}{4}.
\end{equation}
In the case of  $\gamma = 2$, the two linear independent solutions of equation \eqref{indStokes1} 
admit the following representation,
\begin{equation}\label{gamma2}
v^{+}(z) = \Psi_{21}(0, -z6^{1/5}), \quad v^{-}(z) = \Psi_{22}(0, -z6^{1/5}),
\end{equation}
where $\Psi(\lambda, x)$ is the $2\times 2$ matrix solution of the Riemann-Hilbert
problem associated with the first Painlev\'e equation. This Riemann-Hilbert problem
is used in \cite{Ka} for evaluation of the Stokes constant $S_1$. Equations \eqref{gamma2}
allow to use the same Riemann-Hilbert problem to evaluate the Stokes parameter 
associated with the functions $v^{\pm}(z)$. 

The Stokes constant corresponding to the case $\gamma =3/4$ was conjectured in
\cite{GM} on the basis of the relation of equation \eqref{indStokes1} with $\gamma =3/4$
to the symmetric quartic matrix model studied by Br\'ezin--Neuberger \cite{BN}
and Harris--Martinec \cite{HM}. This relation suggests the existence of an alternative
Riemann-Hilbert representation of the first Painlev\'e transcendent which, through the
formulae similar to \eqref{gamma2} would generate solutions of equation \eqref{indStokes1} 
with $\gamma = 3/4$. We expect this alternative Riemann-Hilbert problem can be 
deduced with the help of the relevant double scaling limit of the Lax pair for
the skew-orthogonal polynomials associated with the  symmetric quartic matrix model
\cite{BN}. 

\sectiono{Trans-series} 
\label{sec.trans}

\subsection{Alien derivatives and resurgence equations} An alternative 
method to obtain the asymptotic behavior of the instanton solutions to 
the Painlev\'e I equation is based 
on the ideas of resurgence introduced by \'Ecalle. In this approach, one 
first constructs the {\it trans-series} solution to the differential 
equation and computes its {\it alien derivatives}. 
Using these, one can easily deduce the asymptotics by using contour 
deformation 
arguments. This approach has been used in \cite{CK1, CK2,CG1}, see also 
\cite{GM}, section 4, for an application to 
a first-order differential equation of the Riccati type. In the case of 
Painlev\'e I we deal with a {\it resonant} equation, i.e. there is an integer linear combination 
of its eigenvalues which is null. In this case, the method of 
resurgence is slightly more subtle. In particular, 
the relevant trans-series solution involves logarithms, as we will see. 

As we mentioned in the Introduction, the instanton solutions to Painlev\'e I 
appear as trans-series solutions, see (\ref{uzC}). In order to understand 
their asymptotics we need however 
a more general trans-series solution, which is a formal, 
{\it two}-parameter series of the form 
\be
\label{uzCC}
u(z,C_1, C_2)=\sum_{n,m \ge 0}  C_1^n C_2^m u_{n|m}(z).
\ee
The $u_{n|m}(z)$ have the structure
\be
u_{n|m}(z) =z^{1/2} \re^{-(n-m) A z^{5/4}} \phi_{n|m}(z), 
\ee
where $\phi_{n|m}(z)$ do not contain exponentials $\re^{\pm A z^{5/4}}$. 
When $C_2=0$, $C_1=C$ we recover the trans-series solution (\ref{uzC}), 
therefore
\be
u_{n|0} (z)=u_{n}(z). 
\ee
Notice that the more general trans-series (\ref{uzCC}) is not proper 
(in the sense explained in the Introduction), since for any 
given direction in the complex plane, there is an infinite number of terms 
in (\ref{uzCC}) involving exponentials which grow big 
as $z\rightarrow \infty$ along that direction. Therefore, the trans-series 
(\ref{uzCC}) is more general than those considered in for example \cite{C}. 

If we substitute (\ref{uzCC}) into (\ref{PI}) and collect the coefficients 
of $C_1^nC_2^m$, we find that the $u_{n|m}(z)$ obey a set of coupled, 
inhomogeneous linear ODEs
\be
u''_{n|m} -6 \sum_{k,l\ge0} u_{k|l} u_{n-k| m-l}=0
\ee
generalizing (\ref{coupledC}). In the following it will be convenient to 
introduce the variable
\be
\label{xvar}
x=z^{5/4}.
\ee

We now compute the alien derivatives of these solutions. The alien 
derivative $\Delta_{\omega}$ of a formal power series $\phi(x)$ was introduced 
by \'Ecalle \cite{Ec1,Ec2,Ec3}. To obtain $\Delta_{\omega} \phi$, one 
essentially computes the discontinuity of 
the Borel transform of $\phi(x)$, $\hat \phi(\zeta)$, at the cut in the 
Borel plane starting at 
$\zeta=\omega A$. This leads under suitable assumptions to a series in 
$\IC\{ \zeta\}$, whose inverse Borel transform is the alien derivative 
$\Delta_{\omega} \phi$, 
see also \cite{CNP,SS} for more precise definitions. A crucial property 
is that the pointed alien derivative 
\be
\dot\Delta_\omega =\re^{\omega A x }\Delta_\omega 
\ee
commutes with the standard derivative. This makes possible to relate alien 
derivatives to trans-series solutions. In our case, if we apply the 
pointed alien derivative  
to the Painlev\'e I equation we obtain the linear, second order ODE 
\be
-{1\over 6} {\rd^2 \over \rd z^2} (\dot \Delta_{\omega} u(z,C_1, C_2)) 
+ 2 u(z,C_1, C_2) \dot \Delta_{\omega} u(z,C_1, C_2)=0.
\ee
This equation has two linearly independent solutions:
\be
{\partial u (z,C_1, C_2) \over \partial C_1}, \qquad  
{\partial u (z,C_1, C_2) \over \partial C_2},\ee
and we conclude that 
\be
\dot \Delta_{\omega} u(z,C_1, C_2)=a_{\omega}(C_1, C_2) 
{\partial u (z,C_1, C_2) \over \partial C_1} + b_{\omega}(C_1, C_2) 
{\partial u (z,C_1, C_2) \over \partial C_2}.
\ee
where $a_{\omega}(C_1, C_2)$, $b_{\omega}(C_1, C_2)$ are in principle 
formal power series in $C_1, C_2$ but they are independent of $z$. This 
type of equation, relating the pointed 
alien derivative to trans-series solutions, is called in \'Ecalle's 
theory the {\it bridge equation}; see \cite{CNP,SS} and particularly \cite{GS} for an example of 
a second-order difference equation. In order to understand the 
asymptotics, we are are particularly interested in the cases $\omega=\pm 1$. 
Let us first analyze the case $\omega=1$. We find
\be
\ba
\sum_{n,m\ge 0} \re^{-{(n+1-m)A x}} C_1^n C_2 ^m \Delta_1 \phi_{n|m} 
&= a_1(C_1, C_2) \sum_{n\ge 1, m\ge 0} n C_1^{n-1} C_2^m \re^{-{(n-m)A  x}} 
\phi_{n|m} \\
&+b_1(C_1, C_2)  \sum_{m\ge 1, n\ge 0} m C_1^{n} C_2^{m-1} \re^{-{(n-m)A x}} 
\phi_{n|m}.
\ea
\ee
Let us first look at the term multiplying $\re^{-Ax}$ in both sides. In the 
left hand side this corresponds to $n=m$, and we obtain 
\be
\ba
\sum_{n\ge0} (C_1 C_2)^n \Delta_1 \phi_{n|n}&=a_1(C_1, C_2) 
\sum_{n\ge  0} (n+1) (C_1 C_2)^n  \phi_{n+1|m} \\
&+ b_1(C_1, C_2) \sum_{m\ge  1} m C_1^{m+1} C_2^{m-1}  \phi_{m+1|m}.
\ea
\ee
Since the left hand side is only a function of $C_1 C_2$, we conclude that
\be
a_1(C_1, C_2) =\sum_{k\ge 0} a_{1,k} (C_1 C_2)^k, \qquad b_1(C_1, C_2)
= C_2^2 \sum_{k\ge 0} b_{1,k} (C_1 C_2)^k.
\ee
If we now equate the different powers of $\re^{n  A x}$ and $C_1$, $C_2$, 
we find the equation 
\be
\Delta_1 \phi_{n|m}=\sum_{k=0}^{{\rm min}(n,m)} a_{1,k} (n+1-k) 
\phi_{n+1-k|m-k} +\sum_{k=0}^{{\rm min}(n,m-1)} b_{1,k} (m-1-k) 
\phi_{n-k|m-1-k}.
\ee

The case $\omega=-1$ is very similar. We have now that
\be
a_{-1}(C_1, C_2) =C_1^2 \sum_{k\ge 0} a_{-1,k} (C_1 C_2)^k, 
\qquad b_{-1}(C_1, C_2)= \sum_{k\ge 0} b_{1,k} (C_1 C_2)^k
\ee
and
\be
\Delta_{-1} \phi_{n|m}=\sum_{k=0}^{{\rm min}(n-1,m)} a_{-1,k} (n-1-k) 
\phi_{n-1-k|m-k} +\sum_{k=0}^{{\rm min}(n,m+1)} b_{-1,k} (m+1-k) 
\phi_{n-k|m+1-k}.
\ee

We are particularly interested in the solutions with $m=0$, corresponding to 
the instanton solutions of Painlev\'e I. In this case, the equations for the 
alien derivatives read (we denote $\phi_{k|0}=\phi_k$)
\be
\label{alienone}
\Delta_{1} \phi_k = S_{1}(k+1) \phi_{k+1}, \quad k\ge 0
\ee
where $S_{\pm 1}=a_{\pm1,0}$, and 
\be
\label{alientwo}
\Delta_{-1} \phi_k = S_{-1} (k-1) \phi_{k-1} + \widetilde S_{-1} \phi_{k|1}, 
\qquad k\ge 0,
\ee
where $\widetilde S_{-1}=b_{-1,0}$ and we understand that $\phi_{-1}=0$ in 
the case $k=0$. The three unknown constants $S_{\pm1}$, $\widetilde S_{-1}$ 
are 
the Stokes constants for the Painlev\'e I equation. 

As we will see, a consequence of these equations is that the asymptotics of 
the instanton solution $u_k$ to Painlev\'e I is determined by the solution 
$u_k'$, with $k'=k\pm 1$, and by the 
trans-series solutions $u_{k|1}$, which do not belong to the instanton 
sequence.

\subsection{A study of the $u_{n|1}$ trans-series} 

We will now study the trans-series solutions $u_{n|1}$. It turns out that 
there are {\it three} differerent cases: $n=0$, $n=1$ and $n\ge 2$. We will 
now study them in detail. 

For $n=0$, $u_{0|1}$ satisfies the same linear ODE than 
$u_1$, 
\be
-{1\over 6} u_{0|1}'' + 2u_0 u_{0|1}=0,
\ee
which indeed has two linearly independent solutions: one of them, 
corresponding to $u_1$, 
is exponentially decreasing along the direction ${\rm arg}\, z=0$, 
$|z|\rightarrow \infty$. The solution 
corresponding to $u_{0|1}$ is {\it exponentially increasing} along this 
direction, and it is given by 
\be
u_{0|1}(z)= z^{-1/8} \re^{+ {8 {\sqrt{3}} \over 5} z^{5/4}} \mu_1(z), 
\ee
where $\mu_1(z)$ is a formal power series in $z^{-5/4}$:
\be
\mu_1(z) = \sum_{n\ge 0}  \mu_{n,1} z^{-5n/ 4}
\ee
and we normalize again
\be
\mu_{0,1}=1.
 \ee
The coefficients $ \mu_{n,1}$ satisfy the same recursion than $u_{n,1}$, 
(\ref{oneirec}), with the only difference that we have a minus sign in the 
right hand side. One immediately finds that
\be
\label{relation}
 \mu_{n,1} =(-1)^n u_{n,1}. 
\ee

Let us now consider $u_{1|1}$, which satisfies the linear inhomogeneous ODE
\be
-{1\over 6} u_{1|1}'' + 2u_0 u_{1|1}+ 2 u_1 u_{0|1}=0.
\ee
It is easy to see that $u_{1|1}(z)$ has the following structure
\be
u_{1|1}(z)=z^{-3/4}\sum_{n\ge0} \mu_{n,2} z^{-5n/4}, 
\ee
and in addition $\mu_{2n+1, 2}=0$. The even coefficients satisfy the recursion
\be
\label{mutwo}
\mu_{2n,2}=-\sum_{l=0}^{n-2} \mu_{2l,2} u_{n-l,0} -\sum_{l=0}^{2n}  
(-1)^l u_{l,1} u_{ 2n-l,1} +{25\over 192}(2n-1)^2 \mu_{2(n-1),2}
\ee
and we find, for the very first terms, 
\be
\label{twostrange}
u_{1|1}(z)= -z^{-3/4} \Bigl( 1 +{75\over 512} z^{-5/2} +{300713\over 
1572864}z^{-5} +\cdots\Bigr). 
\ee

For $k\ge 2$, the formal solutions $u_{k|1}$ develop a new feature: 
they contain logarithms. This is due to the resonant character of the 
Painlev\'e I equation. These solutions have the following form: 
\be
\label{logsols}
u_{k|1}(z) ={5\over 4} \log z\, g_{k-1} (z)+f_{k+1} (z), \qquad k\ge 2, 
\ee
where
\be
\label{fg}
\ba
f_k (z)&=z^{1/2 -5 k/8} \re^{(2-k) A z^{5/4}} \mu_{k}(z), 
\quad \mu_{k}(z)=\sum_{n \ge 0} \mu_{n,k}z^{-5n/4}, \\
g_k(z) &=z^{1/2-5k/8} \re^{-k A z^{5/4}} \nu_{k}(z), 
\quad \nu_{k}(z)= \sum_{n\ge0}\nu_{n,k}z^{-5n/4}.
\ea
\ee
The factor $5/4$ in (\ref{logsols}) is introduced for convenience, since 
in the resurgent analysis it will be convenient to use the variable $x$ in 
(\ref{xvar}). 
The functions $f_k, g_k$ appearing in (\ref{fg}) satisfy the coupled system of 
equations
\be
\ba
&-{1\over 6} g_k''+2 u_0 g_k +2\sum_{i=1}^{k-1} u_i 
g_{k-i}=0, \qquad k\ge 1.\\
&-{1\over 6} f_k'' + 2u_0 f_k+ 2 \sum_{i=1}^{k-1} 
u_i  f_{k-i}  + {5\over 24 z^2} g_{k-2}  -{5\over 12 z} 
g_{k-2}'=0,  \qquad k\ge 3.
\ea
\ee
In the second equation we set $f_1(z)=u_{0|1}(z)$, $f_2(z)=u_{1|1}(z)$. 
It is easy to see, from the recursion relation obeyed by the coefficients 
$\nu_{n,k}$, that
 \be
 \label{phiugen}
\nu_{k}(z)=C k u_{k}(z), \qquad k\ge 1
 \ee
where $C$ is a constant given by 
\be
\label{Cvalue}
C={16 \over 5 A}. 
\ee
The value of this constant can be fixed by looking at the equation for 
$u_{2|1}(z)$.

Finally, one can easily find recursion relations for the coefficients 
$\mu_{n,k}$ appearing in $f_{k}(z)$. The cases 
$k=3$ and $k\ge 4$ are slightly different. For $k=3$, one finds 
\be
\label{muthree}
\ba
\mu_{n,3}&={8 \over 25 A (n+1)}\biggl\{ -{25\over 64} (2n+1)^2 \mu_{n-1,3} 
+ 12 \sum_{l=0}^{n-2} u_{(n+1-l)/2,0}\mu_{l,3} \\
& \qquad + 12 \sum_{m=1}^{2} \sum_{l=0}^{n+1} u_{n+1-l,m}\mu_{l,3-m} 
+ {25\over 16} (2n+1) \nu_{n,1} + {25A\over 8} \nu_{n+1,1}\biggr\}, 
\ea
\ee
while for $k\ge 4$ we have
\be
\label{mufour}
\ba
\mu_{n,k}&= \frac{1}{12( k-1)(k-3)} \Bigl\{ 12\sum_{l=0}^{n-3} \mu_{l,k} 
u_{(n-l)/2,0} + 
  12 \sum_{m=1}^{k-1} \sum_{l=0}^{n} \mu_{l,m} u_{n-l,k-m} 
-\frac{25}{64}(2n+k-4)^2 \mu_{n-2,k}   \\
  &
  -\frac{25}{16} A k (k+2n-3) 
\mu_{n-1,k} + {25\over 16} (k+2n-4) \nu_{n-1,k-2} + {25A\over 8}(k-2) 
\nu_{n,k-2}
 \Bigr\}
\ea
\ee

One interesting aspect of 
the doubly-indexed sequences $(u_{n,k})$, $(\mu_{n,k})$ is that we can 
also consider their
dependence on $k$ for fixed $n$, and it turns out that one can find 
closed formulae for these coefficients from the recursion. 
For example, $(u_{0,k})$ is given by
\be
\label{onel}
u_{0,k}=(12)^{1-k} k, \quad k\ge 1,
\ee
while $(u_{1,k})$ is given by
\be
u_{1,k}=-{12^{-k} \over 16 {\sqrt{3}}} \Bigl( 109 \, k^2 -120 
k+ 24 \Bigr), \qquad k \ge 2. 
\ee
These formulae can be also obtained from the results in section 5.2 
of \cite{CC}. Using their results one can see that, as a function of 
$k$,  $12^k u_{n,k}$ is a polynomial in $k$ of degree $n$. 
Finally, we give a general formula for $\mu_{0,k}$ when $k\ge 4$:
\be
\label{mugen}
\mu_{0, k}=12^{-k+1} (k-2)(141 k -402), \qquad k \ge 4. 
\ee
Let us close this section with a problem dealing with the physical interpretation of the full trans-series (\ref{uzCC}). 

\begin{problem}
As we mentioned in the Introduction, the series $u_{n|0}(z)$ can be interpreted in terms of the double-scaling 
limit of instantons in the matrix model, and as an amplitude associated to a ZZ brane in Liouville gravity. What is the interpretation of the 
more general trans-series $u_{n|m}(z)$, with $m>0$, in the context of matrix models and in the context of non-critical strings? Do they 
correspond to new non-perturbative sectors of these theories?
\end{problem}

\sectiono{Asymptotics}
\label{sec.asymptotics}

\subsection{Asymptotics of the multi-instanton solutions}

The asymptotics of the coefficients $u_{n,k}$ can be obtained from 
a well-known application of Cauchy's theorem in the Borel plane (see 
\cite{CK2,CG1,GM} for examples). To do this, we consider the formal power 
series $\phi_k$ appearing in the $k$- instanton solution, which we write in 
terms of the variable $x=z^{5/4}$ as
\be
\phi_k(x)=x^{-\beta k} \sum_{n\ge 0} u_{n,k} x^{-n}
\ee
Here, $\beta=1/2$. The Borel transform of $\phi_k(x)$ is
\be
\hat \phi_k(p)=\sum_{n=0}^{\infty} {u_{n,k} \over \Gamma(\beta k +n)} 
p^{n+\beta k-1}. 
\ee
Therefore, 
\be
{u_{n,k} \over \Gamma(\beta k +n) }=\oint_0{\rd p\over 2\pi \ri} 
{\hat \phi_k (p) \over p^{\beta k+n}}. 
\ee
We now apply the standard deformation contour argument. The Borel transform 
$\hat \phi_k$ has branch cuts at $p=\pm A$, and we can deform the contour to 
encircle these (we also 
pick a vanishing contribution from a circle at infinity). But the 
contribution to the integral is given precisely by the discontinuity 
across the cut, i.e.  
by the alien derivatives that we calculated in (\ref{alienone}) and 
(\ref{alientwo}). We then obtain the asymptotic formula
\be
\label{cdef}
{u_{n,k} \over \Gamma(\beta k +n) }\sim_n {S_1  \over 2 \pi \ri} (k+1) 
\int_{A}^\infty {\rd p} {\hat \phi_{k+1} (p) \over p^{\beta k+n}}
-{S_{-1}  \over 2 \pi \ri} (k-1) \int_{-\infty}^{-A} {\rd p} 
 {\hat \phi_{k-1} (p) \over p^{\beta k+n}}-{\widetilde S_{-1}  
\over 2 \pi \ri}  \int_{-\infty}^{-A} {\rd p} 
 {\hat \phi_{k|1} (p) \over p^{\beta k+n}}.
 \ee
For $k\ge 2$ this formula involves the Borel transform of $\phi_{k|1}(x)$. 
To compute this we have to use
\be
\label{borellog}
\CB \Bigl( \log x \, x^{-\nu} \Bigr)= -{p^{\nu-1} \over \Gamma(\nu)} 
\log p + {\psi(\nu) \over \Gamma(\nu)} p^{\nu-1}. 
\ee
The calculation of the integrals appearing in (\ref{cdef}) is very similar 
to the calculation in \cite{GM}. The only new ingredient is the logarithm 
appearing 
in the Borel transform (\ref{borellog}), which leads to an integral of the 
form 
\be
\int_{0}^{\infty}  \rd \zeta 
\log \zeta {\zeta^{r+\beta (k-1)-1} \over (1+ \zeta)^{k\beta+ n}} 
={ \Gamma(r+\beta (k-1))\Gamma (n-r+\beta)  \over   \Gamma (n+\beta k)} 
\bigl(\psi(n-r+\beta)-\psi ( r+ (k-1) \beta) \bigr).
\ee
Notice that the term involving $\psi ( r+ (k-1) \beta)$ will cancel against 
the contribution coming from the second term in (\ref{borellog}). For $n$ 
large we can use the asymptotic behavior
\be
\label{digas}
 \psi(n-r+\beta)=\log n +{ \beta -r -1/2 \over n} +\CO(1/n^2). 
 \ee

 Before presenting the general result for the asymptotics derived from 
(\ref{cdef}), let us analyze in some detail the case $k=0$, since this 
will fix one of the Stokes constants. 
 For $k=0$ the only contributions to the asymptotics come from $u_1$ and 
$u_{0|1}$, and we obtain 
\be
\ba
u_{n/2,0} & \sim  A^{-n+1/2} {S_1\over 2 \pi\ri}  
\Gamma\bigl(n-1/2 \bigr)\, 
\biggl\{ u_{0,1} + \sum_{l=1}^{\infty} {u_{l,1} A^{l} \over 
\prod_{m=1}^{l} 
(n-1/2 -m)} \biggr\} \\ 
& + 
(-A)^{-n+1/2} {\widetilde S_{-1} \over 2\pi\ri}  \Gamma\bigl(n-1/2 \bigr)\, 
\biggl\{ \mu_{0,1}+ \sum_{l=1}^{\infty} {\mu_{l,1} (-A)^{l} \over 
\prod_{m=1}^{l} 
(n-1/2-m)} \biggr\}.
\ea
\ee
Since $u_{n/2,0}=0$ if $n$ is not even, the right hand side of this 
relation must vanish if $n$ is odd. Using (\ref{relation}) we find that 
this is the case provided that 
\be
\label{tildess}
(-1)^{1\over 2} \widetilde S_{-1} =S_1
\ee
This is exactly as in the ODE studied in 
\cite{CCK}. We then find the result \cite{JK} (see also [GLM]) 
\be
u_{n,0}  \sim   A^{-2n+1/2} {S_1\over  \pi\ri}  
\Gamma\bigl(2n-1/2 \bigr)\, 
\biggl\{1+ \sum_{l=1}^{\infty} {u_{l,1} A^{l} \over 
\prod_{m=1}^{l} 
(2n-1/2 -m)} \biggr\}, \qquad n \rightarrow \infty. 
\ee
Notice that $S_1$ has already been evaluated in \eqref{residue0_integral} and it has the value
\be
S_1 =-\ri {3^{1/4} \over 2 \pi^{1/2}}. 
\ee

Let us now analyze $k=2$. The asympotics of $u_{n,1}$ involves $u_0$ and 
$u_{1|1}$. Using (\ref{tildess}) we can write it as
\be
\label{oneias}
u_{n,1}  \sim A^{-n+1/2} {S_1\over 2 \pi\ri}  \Gamma\bigl(n-1/2 
\bigr)\Bigl\{ 2u_{0,2} + (-1)^n \mu_{0,2}  +  \sum_{l=1}^{\infty} 
{(2 u_{l,2} +(-1)^{n+l}\mu_{l,2}) A^{l} \over \prod_{m=1}^{l} 
(n-1/2 -m)} \biggr\}, 
\ee
and it depends on the parity of $n$. Using the explicit results (\ref{exonei}), 
it is easy to check that the leading asymptotic behavior is as stated by 
Joshi and Kitaev in \cite[Prop.16]{JK}. The formula (\ref{oneias})
gives in addition the all-orders expansion of $u_{n,1}$ as an asymptotic 
series in $1/n$. 

We can now write a general formula for the asymptotics of $u_{n,k}$ when 
$k\ge2$. Using again (\ref{tildess}), and absorbing a factor $(-1)^{-1/2}$ 
in $S_{-1}$, 
we find
\be
\label{asunk}
\ba
u_{n,k} & \sim  A^{-n+\beta} {S_1\over 2 \pi\ri}  
\Gamma\bigl(n-\beta \bigr)\, 
\biggl\{ (k+1) u_{0,k+1} +(-1)^n \mu_{0,k+1}  
+ \sum_{l=1}^{\infty} {((k+1)u_{l,k+1} +(-1)^{n+l} 
\mu_{l,k+1}) A^{l} \over \prod_{m=1}^{l} 
(n-\beta -m)} \biggr\} \\ 
& + 
(-1)^n (k-1) A^{-n-\beta} {S_{-1} \over 2\pi\ri}  
\Gamma\bigl(n+\beta \bigr)\, 
\biggl\{ u_{0,k-1}+ \sum_{l=1}^{\infty} {u_{l,k-1} (-A)^{l} \over 
\prod_{m=1}^{l}
(n+\beta-m)} \biggr\}\\
&- (-1)^n A^{-n-\beta} {S_{1} \over 2\pi\ri}  
\Gamma\bigl(n+\beta \bigr)\, (\log n -\log A) \biggl\{ \nu_{0,k-1}
+ \sum_{l=1}^{\infty} {\nu_{l,k-1} (-A)^{l} \over \prod_{m=1}^{l}
(n+\beta-m)} \biggr\}\\
 &-(-1)^n A^{-n-\beta} {S_{1} \over 2\pi\ri}  \Gamma\bigl(n+\beta 
\bigr)\,  \biggl\{ (\psi(n+\beta) -\log n) \nu_{0,k-1} 
+ \sum_{l=1}^{\infty} {\psi(n+\beta-l) 
-\log n\over \prod_{m=1}^{l}
(n+\beta-m)} \nu_{l,k-1} (-A)^{l} \biggr\}
\ea
\ee
for $k\ge 2$, and we recall that $\beta=1/2$. As compared to the 
asymptotics for $k=0,1$, the asymptotics for $k\ge2$ involves 
logarithmic terms. In fact, 
the dominant term in the asymptotics is precisely the $\log n$ term. In the 
last line we use the asymptotics (\ref{digas}) for the $\psi$ function. 
Finally, notice from (\ref{phiugen}) that one has the relation 
(\ref{eq.phirec}) for the 
coefficients $\nu_{n,k}$. 

\subsection{Asymptotics of the multi-instanton solutions: numerical 
evidence}

We will now perform numerical tests of the predicted asymptotic behaviour 
(\ref{asunk}) for the instanton series $u_{n,k}$. The standard  technique 
to do that is the method of Richardson extrapolation. This method goes as follows. Let us assume that a sequence 
$s_n$ has the asymptotics
\be
\label{sequ}
s_n \sim_n \sum_{k=0}^{\infty} {a_k \over n^k}
\ee
for $n$ large. Its 
{\it $N$-th Richardson transformation} $s^{(N)}_n$ can be defined recursively by 
\be
\label{richardson}
\ba
s^{(0)}_n&=s_n, \\
s^{(N)}_n&=s_{n+1}^{(N-1)} + {n \over N}(s_{n+1}^{(N-1)} -s_n^{(N-1)}), 
\quad N\ge 1.
\ea
\ee
The effect of this transformation is to remove subleading tails in 
(\ref{sequ}), and
\be
s^{(N)}_n \sim a_0 + \CO\Bigl({1\over n^{N+1}}\Bigr). 
\ee
 The values $s^{(N)}_n$ give numerical approximations to $a_0$, 
and these 
approximations become better as $N$, $n$ increase. Once a numerical 
approximation to $a_0$ has been obtained, the value of $a_1$ can be estimated 
by considering the sequence $n(s_n-a_0)$, and so on. 

The method of Richardson extrapolation can be applied {\it verbatim} to 
verify the asymptotics of $u_{n,1}$ written down in (\ref{oneias}). 
Let us illustrate this 
with a nontrivial example. According to (\ref{oneias}), the sequence
\be
\label{nll}
s_n=(2n-5/2)\biggr\{(2n-3/2) \Bigl( {u_{2n,1} A^{2n-1/2} \over 
S_1/(2\pi \ri)\Gamma\bigl(2n-1/2 \bigr)} +{2\over 3}\Bigr) -2 u_{1,2} \biggl\}
\ee
is of the form (\ref{sequ}) and asymptotes, as $n \rightarrow \infty$, 
the value 
\be\label{nllcoeff}
(2 u_{2,2} +\mu_{2,2})A^{2} =-{55\over 96} \approx -0.5729166666666666...
\ee

\begin{figure}[htpb]
$$
\psdraw{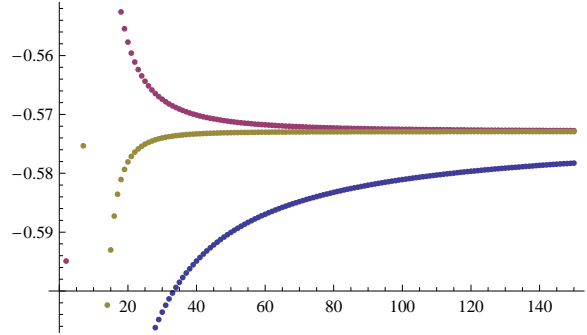}{3in}
$$
\caption{A plot of the sequence (\ref{nll}) and its Richardson transformations. In this and subsequent plots, the horizontal axis represents the integer $n$, and the vertical 
axis the values of the sequence.}
\label{teststrange}
\end{figure}

In \figref{teststrange} we show a plot of the sequence (\ref{nll}) and its 
first two Richardson transformations for $n$ up to $200$,  which converges to the 
expected value (\ref{nllcoeff}). 
Taking $n=250$ and ten Richardson transformations gives the numerical 
approximation
\be
s_{200}^{(5)}=-0.572916666666666667...
\ee
Notice that this test already verifies that the general trans-series solutions 
$u_{n|m}$ appearing in (\ref{uzCC}) are the relevant objects to understand 
the asymptotics, since $k=1$ involves the trans-series $u_{1|1}$.

Let us now study the asymptotic behavior of the instanton sequences 
$u_{n,k}$ with $k\ge 2$. The main novelty here is the presence of 
logarithms in the asymptotics, and 
 this does not fit {\it a priori} into the standard framework of Richardson 
transformations. However, one can transform the sequence and put it in a form 
which is amenable 
 to an analysis with standard Richardson transformations, as pointed out in 
\cite{zj}. Let us assume that we have a sequence $\ell_m$ with the asymptotics
\be
\label{logasym}
\ell_m \sim 
\log m \, s_{m}  + t_{m}, \qquad s_m =\sum_{k\ge 0} {a_k \over m^k}, 
\quad t_m =\sum_{k\ge 0} {b_k \over m^k}
\ee
This type of asymptotc behavior appears in instanton corrections in Quantum 
Mechanics. The leading behavior of this sequence is determined by the 
coefficient $a_0$, and we would like to find a method 
to extract it numerically. To do this we consider the sequence
\be
\tilde \ell_m =m(\ell_{m+1} -\ell_{m}),  
\ee
which has the asymptotics
\be
\tilde \ell_m \sim 
\log m \, \tilde s_{m}  + \tilde t_{m}, \qquad \tilde s_m 
=\sum_{k\ge 1} {\tilde a_k \over m^k}, \quad \tilde t_m 
=a_0+\sum_{k\ge 1} {\tilde b_k \over m^k},
\ee
It is now easy to see that, if we apply the Richardson transformation 
(\ref{richardson}) {\it twice} to this sequence, we remove both the 
tails in $1/n^p$ and the tails in $\log n/n^p$. This then 
allows a precise determination of the leading term $a_0$. Once this has 
been determined, we can extract the other coefficients in (\ref{logasym}) 
by subtracting from the original sequence the parts of the asymptotics 
which are under control. 

Let us now apply this idea to the sequence of instantons of Painlev\'e I. 
The first step is to consider the auxiliary sequence 
\be
a_{n,k} ={A^{n+\beta} u_{n,k} \over \Gamma(n+\beta)}
\ee
whose leading asymptotics is 
\be
\ba
\label{ankas}
(-1)^n a_{n,k}  &\sim  -  {S_{1} \over 2\pi\ri} \biggl\{ \nu_{0,k-1} 
-  {A \nu_{1,k-1} \over n}+\CO(1/n^2)  \biggr\}\, \log n \\
&+  {S_{1} \over 2\pi\ri} \log A \,  \nu_{0,k-1} + {S_{-1} 
\over 2\pi\ri} (k-1)  u_{0,k-1} \\
&+ A \biggl\{  {S_1\over 2\pi \ri}  \Bigl( (k+1) (-1)^n u_{0,k+1} 
+ \mu_{0,k+1}) - \log A \, \nu_{1,k-1}\Bigr)- {S_{-1} \over 2\pi\ri} 
(k-1) u_{1,k-1} \biggr\} {1\over n} \\
& +\CO(1/n^2)
\ea
\ee
Again, it depends on the parity of $n$, i.e. wether $n=2m$ or $n=2m+1$. 
Let us denote
\be
 a^{(e)}_{m,k}=a_{2m,k}, \qquad  a^{(o)}_{m,k}=a_{2m+1,k}.
 \ee
In both cases their asymptotics is of the form (\ref{logasym}) and we can 
use the method of \cite{zj} to analyze the sequence numerically. 

For simplicity, we will illustrate the asymptotics by focusing on the 
sequence with $n=2m$ even, and we will test the four leading terms 
displayed in (\ref{ankas}), i.e. the (leading) term 
in $\log \, n$, and the terms in $\log \, n/n$, constant, and $1/n$. In 
terms of the general structure written down in (\ref{logasym}), we will 
test the values of $a_0, a_1$ and $b_0, b_1$. 
These coefficients involve all the trans-series solutions appearing in the 
resurgent analysis. 

\begin{figure}[htpb]
$$
\psdraw{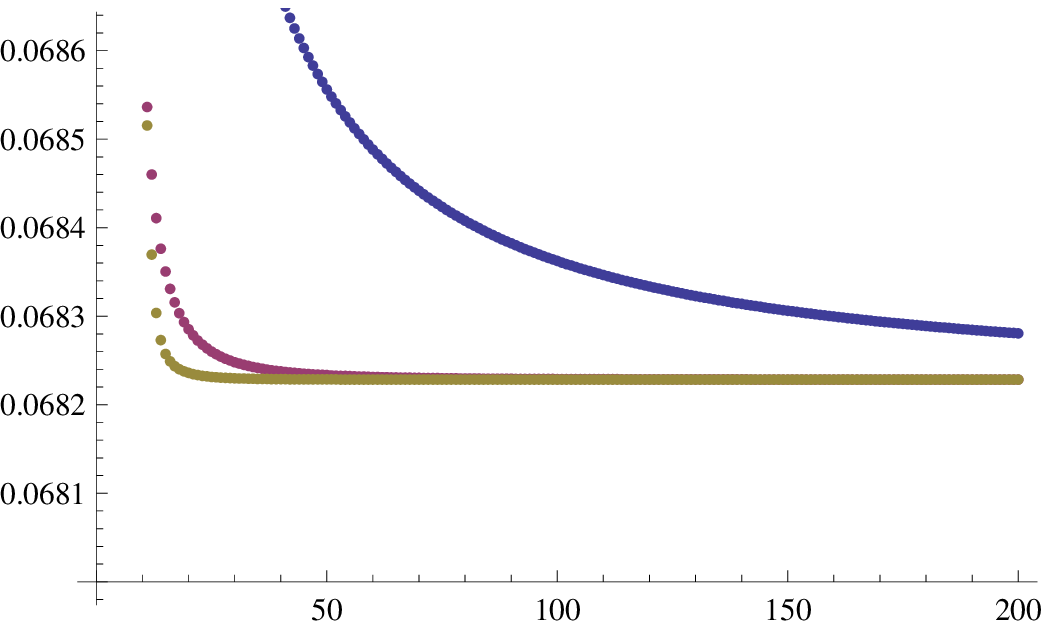}{2.5in} \qquad \psdraw{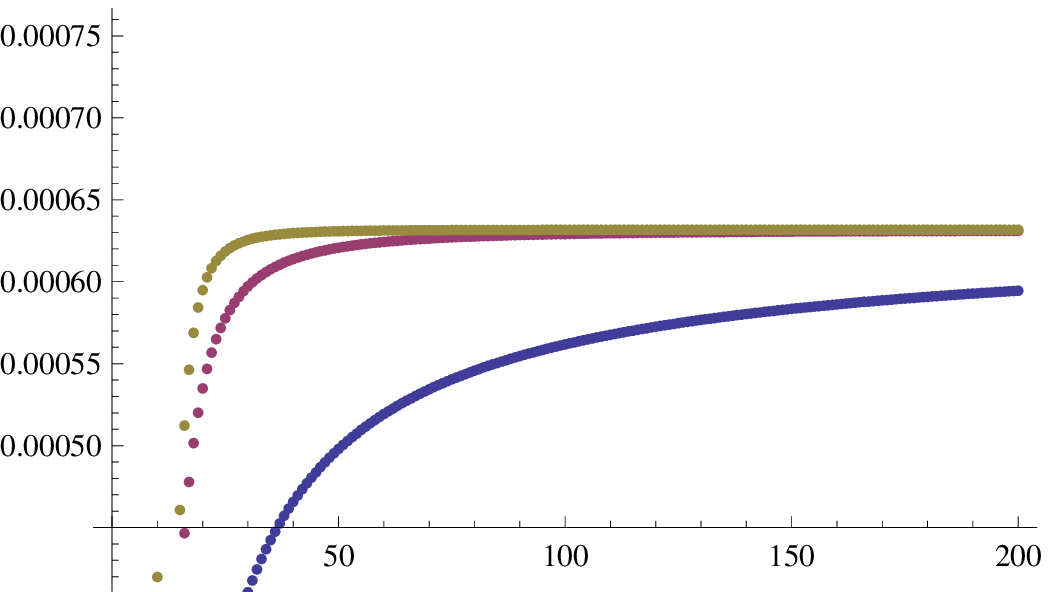}{2.5in}
$$
\caption{A plot of the sequence (\ref{ll}) and its Richardson transformations, 
for $k=2$ (left) and $k=5$ (right).}
\label{twoll}
\end{figure}

We will first test the {\it leading term} of order $\log n$. Following 
\cite{zj}, we first construct the 
sequence 
\be
\label{ll}
\tilde a^{(e)}_{m,k}=s(a_{2m+2,k}-a_{2m,k}) \sim  
-{S_{1} \over \pi\ri {\sqrt{3}}}  (k-1)^2 (12)^{2-k}  
+\CO(1/m, \log m/m)
\ee
where we have used the explicit value for $\nu_{0,k-1}$ derived from 
(\ref{eq.phirec}) and (\ref{onel}). To remove tails, we perform Richardson 
transformations in the sequence 
(each transform is performed twice to remove both types of tails, as 
explained above). In \figref{twoll} we plot it for $k=2,5$, together 
with the second and the fourth Richardson 
transformations, and 
for $m=200$. The convergence towards
\be
 -{S_{1} \over \pi\ri {\sqrt{3}}} \approx 0.068228352037086...
\ee
for $k=2$, and towards 
\be
 -{S_{1} \over 108\pi\ri  {\sqrt{3}}}  \approx 0.00063174400034...
\ee
for $k=5$, is manifest in the figures. More precisely, for $m=250$ and with 
$N=10$ Richardson transformations, we obtain the numerical approximations
\be
\ba
\tilde a_{250,2}^{(e), (10)}&=0.068228352037087...,\\
\tilde a_{250,5}^{(e), (10)}&= 0.00063174400031...
\ea
\ee

\begin{figure}[htpb]
$$
\psdraw{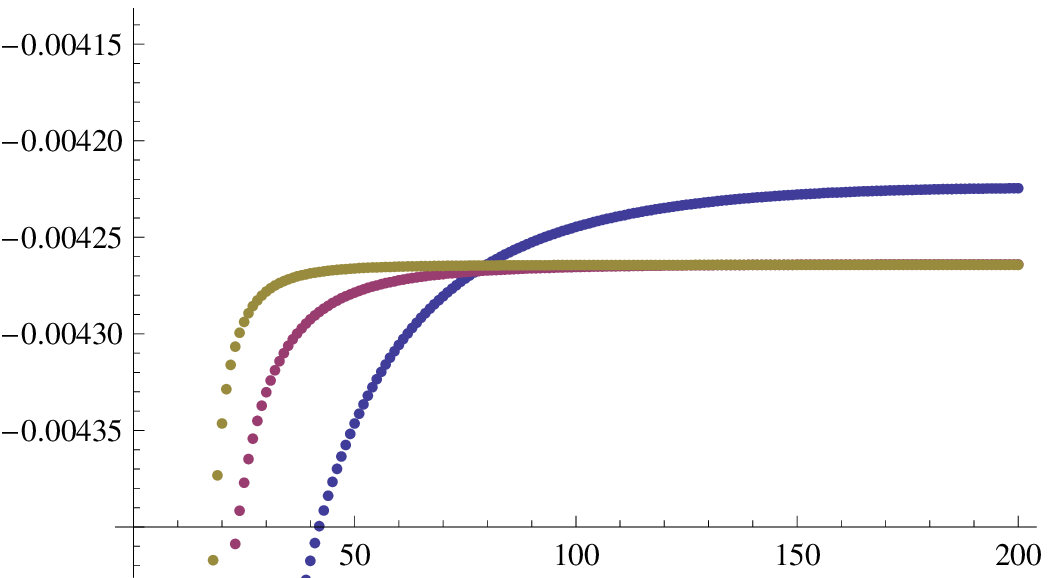}{2.5in} \qquad \psdraw{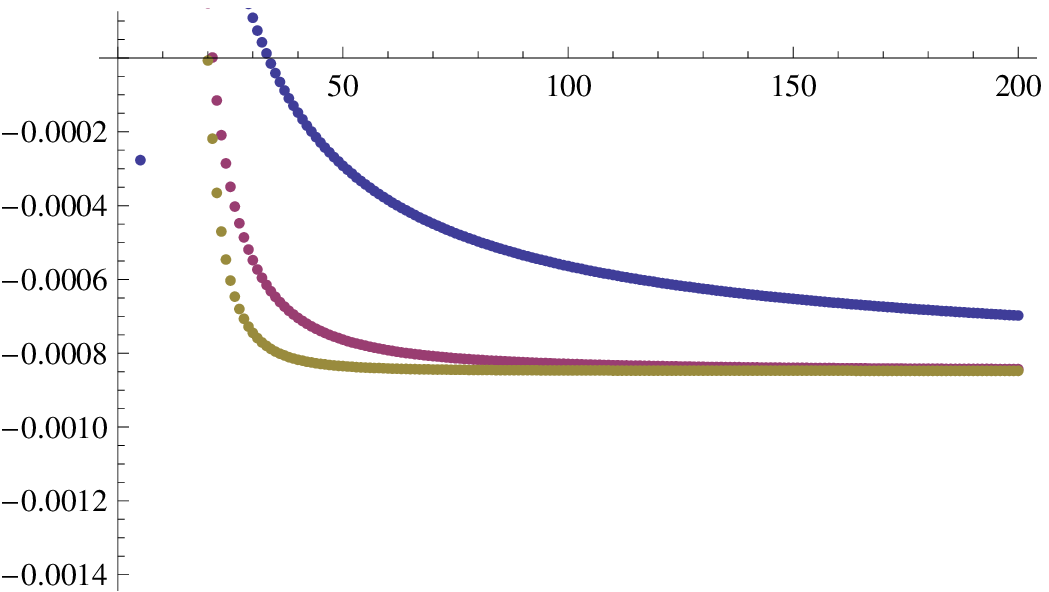}{2.5in}
$$
\caption{A plot of the sequence $\tilde b_{m,k}$, obtained from 
(\ref{nllog}), together with its Richardson transformations, 
for $k=2$ (left) and $k=5$ (right).}
\label{nllogpic}
\end{figure}

We next test numerically the coefficient of $\log s/s$. We consider the 
sequence 
\be
\label{nllog}
b_{m,k}=m\Bigl(\tilde a^{(e)}_{m,k}+ {S_{1} \over 
\pi\ri{\sqrt{3}}} (k-1)^2 (12)^{2-k}  \Bigr)
\ee
whose leading asymptotics is of the form 
\be
-{4 S_{1} \over 5\pi\ri} (k-1) u_{1,k-1} \log m +\cdots
\ee
so we can apply the above procedure. The sequence $\tilde b_{m,k}$ for 
$k=2,5$, and up to $m=200$, 
together with its second and fourth Richardson transformations, is displayed 
in \figref{nllogpic}. 
In both cases we have convergence to the predicted values
\be
k=2:\, \, -0.00426427200231..., \qquad k=5:\, \, -0.000847589867... 
\ee
For $m=250$ and with $N=10$ Richardson transformations, we obtain the numerical 
approximations
\be
\ba
\tilde b_{250,2}^{(10)}&=-0.00426427200235...,\\
\tilde b_{250,5}^{(10)}&=-0.000847589866...
\ea
\ee

We can now study the constant term of the asymptotics, which also makes 
possible to obtain a numerical determination of the 
additional Stokes parameter $S_{-1}$. We consider the sequence
\be
\label{bseq}
c_{m,k}=a_{2m,k} + {S_{1} \over 2\pi\ri} \biggl\{ 
{2\over {\sqrt{3}}}(k-1)^2 (12)^{2-k} 
- {A \nu_{1,k-1} \over 2m} \biggr\}\, 
\log (2m) 
\ee
where we have subtracted the leading logs. According to the predictions of 
resurgence, as $m\rightarrow 
\infty$ this sequence asymptotes to 
\be
\label{cas}
{S_{-1}' \over 2 \pi \ri} (k-1)^2(12)^{2-k}
\ee
where
\be
S_{-1}'= S_{-1}+   {2 \log A \over {\sqrt{3}}}S_{1} .
\ee
The numerical analysis in the case of $k=2$ gives a numerical 
determination of the unknown 
Stokes constant $S'_{-1}$ (hence, of $S_{-1}$), and we find, numerically, 
\be
\label{newstokesval}
 {S_{-1}' \over 2 \pi \ri} \approx 0.31873285573864121...
 \ee
and consequently
\be
\label{S-1}
{S_{-1} \over 2 \pi \ri} \approx 0.3882786818052856841...
\ee
Using this value, we can verify the asymptotic behavior (\ref{cas}) for 
the sequence (\ref{bseq}) for higher values of $k$. 

\begin{figure}[htpb]
$$
\psdraw{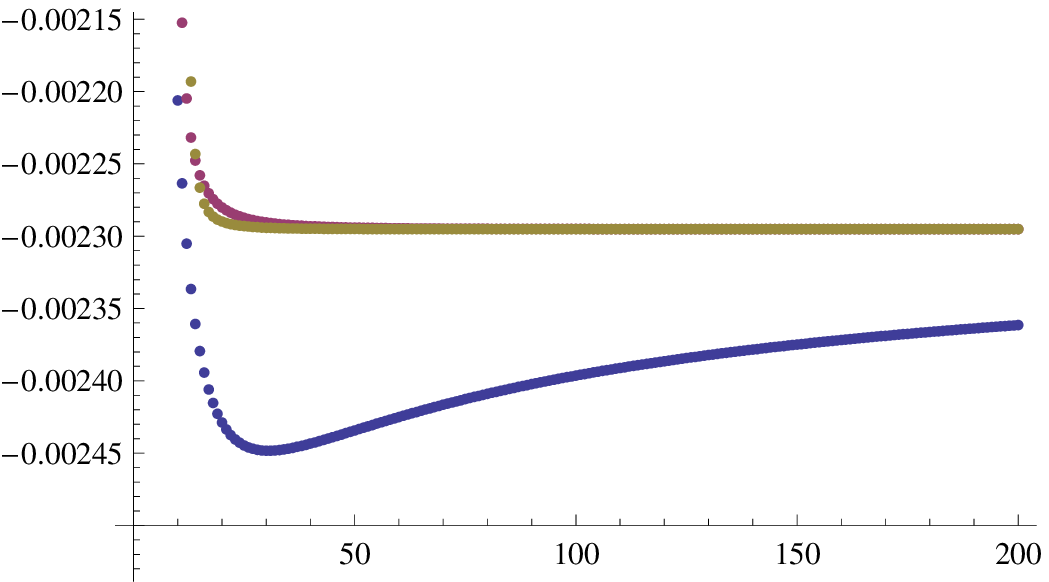}{2.5in} \qquad \psdraw{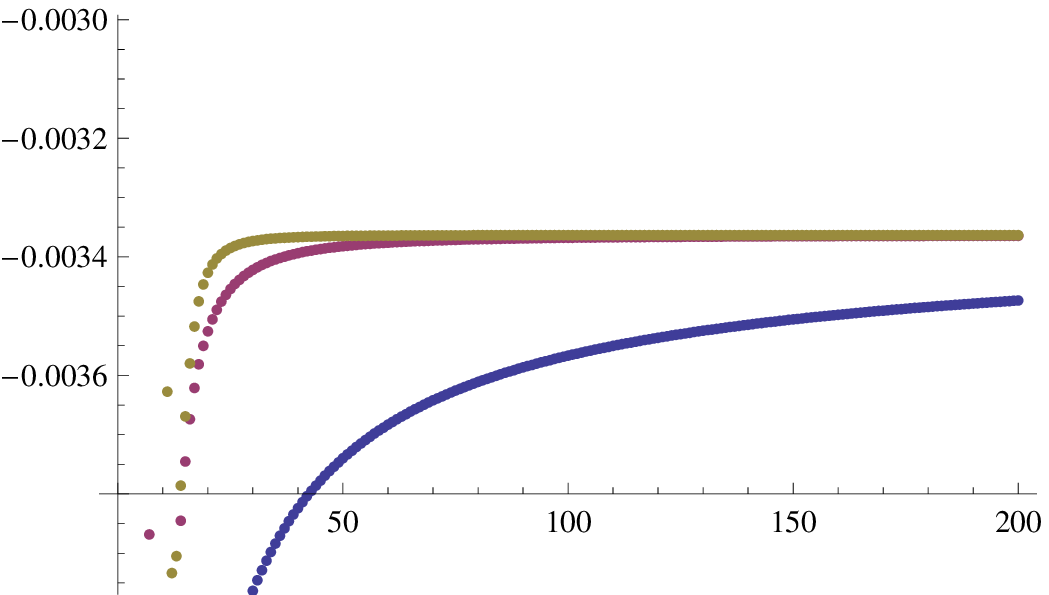}{2.5in}
$$
\caption{A plot of the sequence (\ref{bess}) and its Richardson transformations, 
for $k=2$ (left) and $k=5$ (right).}
\label{fourl}
\end{figure}

Finally, we consider the coefficient of $1/m$. To do this, we construct the 
sequence
\be
\label{bess}
d_{m,k}=m^2(b_{m+1,k}-b_{m,k})
\ee
which should asymptote
\be
-{A \over 2} \biggl\{  {S_1\over 2\pi \ri}  \Bigl( (k+1) u_{0,k+1} 
+ \mu_{0,k+1})  -{S'_{-1} \over 2\pi\ri}  (k-1) u_{1,k-1} \biggr\}. 
\ee
This involves the tran-series coefficients $\mu_{0,k+1}$, and we use the 
numerical determination of $S'_{-1}$ obtained above. In \figref{fourl} we 
show the sequence 
(\ref{bess}) for $k=2$ and $k=5$, up to $m=200$, 
together with its second and fourth Richardson transformations. They clearly 
asymptote the predictions of 
resurgence, 
\be
k=2:\, \, -0.002295145874084..., \qquad k=5: \, \, -0.0033633587118...
\ee
For $m=250$ and with $N=10$ Richardson transformations, we obtain the numerical 
approximations
\be
\ba
d_{250,2}^{(10)}&=-0.002295145874083...,\\
d_{250,5}^{(10)}&=-0.0033633587119...
\ea
\ee

We believe these numerical tests confirm in a very clear way the predictions 
from the resurgent analysis for $k\ge 2$.

\bibliographystyle{hamsalpha}\bibliography{biblio}

\def\cprime{$'$}
\providecommand{\bysame}{\leavevmode\hbox to3em{\hrulefill}\thinspace}
\providecommand{\href}[2]{#2}
\providecommand{\eprint}{\begingroup \urlstyle{rm}\Url}
\begin{thebibliography}{DFGZJ95}

\bibitem[AKK03]{akk}
Sergei~Yu. Alexandrov, Vladimir~A. Kazakov, and David Kutasov,
  \emph{Non-perturbative effects in matrix models and {D}-branes}, J. High
  Energy Phys. (2003), no.~9, 057, 25 pp. (electronic).

\bibitem[BN91]{BN}
{\'E}.~Br{\'e}zin and H.~Neuberger, \emph{Multicritical points of unoriented
  random surfaces}, Nuclear Phys. B \textbf{350} (1991), no.~3, 513--553.

\bibitem[CC01]{CC}
O.~Costin and R.~D. Costin, \emph{On the formation of singularities of
  solutions of nonlinear differential systems in antistokes directions},
  Invent. Math. \textbf{145} (2001), no.~3, 425--485.

\bibitem[CCK04]{CCK}
Ovidiu Costin, Rodica~D. Costin, and Matthew Kohut, \emph{Rigorous bounds of
  {S}tokes constants for some nonlinear ordinary differential equations at
  rank-one irregular singularities}, Proc. R. Soc. Lond. Ser. A Math. Phys.
  Eng. Sci. \textbf{460} (2004), no.~2052, 3631--3641.

\bibitem[CG11]{CG1}
Ovidiu Costin and Stavros Garoufalidis, \emph{Resurgence of the
  kontsevich-zagier power series}, Ann. Inst. Fourier (2011).

\bibitem[CK96]{CK1}
Ovidiu Costin and Martin~D. Kruskal, \emph{Optimal uniform estimates and
  rigorous asymptotics beyond all orders for a class of ordinary differential
  equations}, Proc. Roy. Soc. London Ser. A \textbf{452} (1996), no.~1948,
  1057--1085.

\bibitem[CK99]{CK2}
\bysame, \emph{On optimal truncation of divergent series solutions of nonlinear
  differential systems}, R. Soc. Lond. Proc. Ser. A Math. Phys. Eng. Sci.
  \textbf{455} (1999), no.~1985, 1931--1956.

\bibitem[CNP93]{CNP}
B.~Candelpergher, J.-C. Nosmas, and F.~Pham, \emph{Approche de la
  r\'esurgence}, Actualit\'es Math\'ematiques. [Current Mathematical Topics],
  Hermann, Paris, 1993.

\bibitem[Cos98]{C}
Ovidiu Costin, \emph{On {B}orel summation and {S}tokes phenomena for rank-{$1$}
  nonlinear systems of ordinary differential equations}, Duke Math. J.
  \textbf{93} (1998), no.~2, 289--344.

\bibitem[Dav93]{david}
Fran{\c{c}}ois David, \emph{Nonperturbative effects in matrix models and vacua
  of two-dimensional gravity}, Phys. Lett. B \textbf{302} (1993), no.~4,
  403--410.

\bibitem[DFGZJ95]{DGZ}
P.~Di~Francesco, P.~Ginsparg, and J.~Zinn-Justin, \emph{{$2$}{D} gravity and
  random matrices}, Phys. Rep. \textbf{254} (1995), no.~1-2, 133.

\bibitem[{\'E}ca81a]{Ec1}
Jean {\'E}calle, \emph{Les fonctions r\'esurgentes. {T}ome {I}}, Publications
  Math\'ematiques d'Orsay 81 [Mathematical Publications of Orsay 81], vol.~5,
  Universit\'e de Paris-Sud D\'epartement de Math\'ematique, Orsay, 1981, Les
  alg{\`e}bres de fonctions r{\'e}surgentes. [The algebras of resurgent
  functions], With an English foreword.

\bibitem[{\'E}ca81b]{Ec2}
\bysame, \emph{Les fonctions r\'esurgentes. {T}ome {II}}, Publications
  Math\'ematiques d'Orsay 81 [Mathematical Publications of Orsay 81], vol.~6,
  Universit\'e de Paris-Sud D\'epartement de Math\'ematique, Orsay, 1981, Les
  fonctions r{\'e}surgentes appliqu{\'e}es {\`a} l'it{\'e}ration. [Resurgent
  functions applied to iteration].

\bibitem[{\'E}ca85]{Ec3}
\bysame, \emph{Les fonctions r\'esurgentes. {T}ome {III}}, Publications
  Math\'ematiques d'Orsay [Mathematical Publications of Orsay], vol.~85,
  Universit\'e de Paris-Sud, D\'epartement de Math\'ematiques, Orsay, 1985,
  L'{\'e}quation du pont et la classification analytique des objects locaux.
  [The bridge equation and analytic classification of local objects].

\bibitem[FIKN06]{FIKN}
Athanassios~S. Fokas, Alexander~R. Its, Andrei~A. Kapaev, and Victor~Yu.
  Novokshenov, \emph{Painlev\'e transcendents}, Mathematical Surveys and
  Monographs, vol. 128, American Mathematical Society, Providence, RI, 2006,
  The Riemann-Hilbert approach.

\bibitem[GLM08]{GLM}
Stavros Garoufalidis, Thang T.~Q. L{\^e}, and Marcos Mari{\~n}o,
  \emph{Analyticity of the free energy of a closed 3-manifold}, SIGMA Symmetry
  Integrability Geom. Methods Appl. \textbf{4} (2008), Paper 080, 20.

\bibitem[GM10]{GM}
Stavros Garoufalidis and Marcos Mari{\~n}o, \emph{Universality and asymptotics
  of graph counting problems in non-orientable surfaces}, J. Combin. Theory
  Ser. A \textbf{117} (2010), no.~6, 715--740.

\bibitem[GS01]{GS}
V.~Gelfreich and D.~Sauzin, \emph{Borel summation and splitting of separatrices
  for the {H}\'enon map}, Ann. Inst. Fourier (Grenoble) \textbf{51} (2001),
  no.~2, 513--567.

\bibitem[HM90]{HM}
Geoffrey~R. Harris and Emil~J. Martinec, \emph{Unoriented strings and matrix
  ensembles}, Phys. Lett. B \textbf{245} (1990), no.~3-4, 384--392.

\bibitem[JK01]{JK}
N.~Joshi and A.~V. Kitaev, \emph{On {B}outroux's tritronqu\'ee solutions of the
  first {P}ainlev\'e equation}, Stud. Appl. Math. \textbf{107} (2001), no.~3,
  253--291.

\bibitem[Kap88]{Ka1}
A.~A. Kapaev, \emph{Asymptotic behavior of the solutions of the {P}ainlev\'e
  equation of the first kind}, Differentsial\cprime nye Uravneniya \textbf{24}
  (1988), no.~10, 1684--1695, 1835.

\bibitem[Kap04]{Ka}
\bysame, \emph{Quasi-linear stokes phenomenon for the {P}ainlev\'e first
  equation}, J. Phys. A \textbf{37} (2004), no.~46, 11149--11167.

\bibitem[Mar]{martinec}
Emil~J. Martinec, \emph{The annular report on non-critical string theory},
  \eprint{arXiv:hep-th/0305148}, Preprint 2003.

\bibitem[MSW08]{msw}
Marcos Mari{\~n}o, Ricardo Schiappa, and Marlene Weiss, \emph{Nonperturbative
  effects and the large-order behavior of matrix models and topological
  strings}, Commun. Number Theory Phys. \textbf{2} (2008), no.~2, 349--419.

\bibitem[Olv97]{O}
Frank W.~J. Olver, \emph{Asymptotics and special functions}, AKP Classics, A K
  Peters Ltd., Wellesley, MA, 1997, Reprint of the 1974 original [Academic
  Press, New York; MR0435697 (55 \#8655)].

\bibitem[SS03]{SS}
Tere~M. Seara and David Sauzin, \emph{Borel summation and the theory of
  resurgence}, Butl. Soc. Catalana Mat. \textbf{18} (2003), no.~1, 131--153.

\bibitem[Tak00]{Tak}
Yoshitsugu Takei, \emph{An explicit description of the connection formula for
  the first {P}ainlev\'e equation}, Toward the exact {WKB} analysis of
  differential equations, linear or non-linear ({K}yoto, 1998), Kyoto Univ.
  Press, Kyoto, 2000, pp.~204, 271--296.

\bibitem[ZJ81]{zj}
J.~Zinn-Justin, \emph{Expansion around instantons in quantum mechanics}, J.
  Math. Phys. \textbf{22} (1981), no.~3, 511--520.

\end{thebibliography}
\end{document}